\documentclass[12pt]{amsart}
\usepackage{amssymb}
\usepackage{amscd}
\usepackage{amsmath}
\usepackage{amsthm}
\def\poh{{\rho_{H(R/\overline L)}}}
\def\ph{{\rho_{\H}}}
\def\Tor{{\rm Tor}}
\def\F{{\mathcal F}}
\def\K{{\mathcal K}}
\def\E{{\mathcal E}}
\def\N{{\mathbb N}}
\def\Y{{\underline{\boldsymbol y}}}
\def\X{{\underline{\boldsymbol x}}}
\def\a{{\{a_1, \dots , a_n\}}}
\def\ap{{\{a'_1, \dots , a'_n\}}}
\def\ax{{\{x_1^{a_1}, \dots , x_n^{a_n}\}}}

\def\H{{\mathcal H}}
\newcommand{\poly}[1]{k[x_1,\dots ,x_{#1}]}

\def\deg{{\rm deg}}
\def\A{{\mathbb A}}
\def\soc{{\rm soc}}

\def\B{{\mathbb B}}
\def\b{{\{b_1, \dots , b_n\}}}
\def\LP{{\mathcal LP}^{\H}_{\A}}

\def\lppb{{L_{\H,\A}}}
\def\LA{{LPP(\A)}}
\newcommand{\LL}[1]{{\rm lpp}_\le{(#1)}}

\newcommand{\m}[1]{{\langle #1 \rangle_{\A}}}
\newcommand{\T}{\mathcal{T}}
\newcommand{\mS}{\mathcal{S}}
\newcommand{\bb}[1]{\mathbb{#1}}

\newcommand{\pts}{\mathbb{X}}
\newcommand{\proj}{\mathbb{P}}

\newtheorem{theorem}{Theorem}

\newtheorem{defprop}{Definition/Proposition}
\newtheorem{corollary}{Corollary}
\newtheorem{proposition}{Proposition}
\newtheorem{conjecture}{Conjecture}
\newtheorem{lemma}{Lemma}

\theoremstyle{definition}
\newtheorem{definition}{Definition}
\newtheorem{remark}{Remark}
\newtheorem{example}{Example}

\title{The Residuals of Lex Plus Powers Ideals and the
  Eisenbud-Green-Harris Conjecture}
\author{Benjamin P. Richert\hspace{.5cm}\\ \\
Sindi Sabourin}
\begin{document}
\maketitle
\baselineskip .3in

\begin{abstract} 
  The $n$-type vectors introduced by Geramita, Harima and Shin are in
  1-1 correspondence with the Hilbert functions Artinian of lex
  ideals.  Letting $\mathbb{A} =\{ a_1,\ldots ,a_n\}$ define the
  degrees of a regular sequence, we construct
  $\LL{\mathbb{A}}$-vectors which are in 1-1 correspondence with the
  Hilbert functions of certain lex plus powers ideals (depending on
  $\mathbb{A}$).  This construction enables us to show that the
  residual of a lex plus powers ideal in an appropriate regular
  sequence is again a lex plus powers ideal.  We then use this result
  to show that the Eisenbud-Green-Harris conjecture is equivalent to
  showing that lex plus powers ideals have the largest {\it last}
  graded Betti numbers (it is well-known that the
  Eisenbud-Green-Harris conjecture is equivalent to showing that lex
  plus powers ideals have the largest first graded Betti numbers).
\end{abstract}

\section{Introduction}

Hilbert functions, in general, have been extensively studied.  Let
$R=k[x_1,\ldots ,x_n]$, where each $x_i$ has degree 1.  Then F.S.
Macaulay (\cite{Macaulay}) characterized those sequences (called
$O$-sequences) which occur as the Hilbert function of any $k$-algebra
of the form $R/I$, where $I$ is a homogeneous ideal.  He showed that a
sequence $S=\{ c_i\}_{i\geq 0}$ is such a Hilbert function if and only
if $c_{i+1}\leq c_i^{\langle i\rangle}$, where $-^{\langle i\rangle}$,
known as Macaulay's function, is expressed in terms of the
$i$-binomial expansion of an integer.  In proving his result, Macaulay
shows that lex ideals have the largest first graded Betti numbers
among all ideals having a fixed Hilbert function.  Bigatti
(\cite{Bigatti}) and Hulett (\cite{Hulett}) have independently
generalized this by showing that, over fields of characteristic 0, lex
ideals have the largest graded Betti numbers (not just the largest
first graded Betti numbers) among all ideals having a fixed Hilbert
function.  Pardue (\cite{Pardue}) generalized this to fields of
arbitrary characteristic.

At about the same time that Bigatti and Hulett proved their result,
Eisenbud, Green and Harris together conjectured that a generalization
in a different direction of Macaulay's result should be true.  Instead
of restricting their attention to lex ideals, they look at ideals
which, modulo appropriate powers of the variables, are lex ideals.
These ideals have become known as lex plus powers ideals; letting
$\mathbb{A}=\a$ be a list of positive integers with $a_1\leq \ldots
\leq a_n$, an ideal $L$ containing $x_1^{a_1}, \ldots , x_n^{a_n}$ as
minimal generators is an $\mathbb{A}$-lex plus powers ideal if
$\overline{L}$ is a lex ideal in $R/\langle x_1^{a_1}, \ldots
,x_n^{a_n}\rangle$.  The conjecture states that as long as there is an
$\mathbb{A}$-lex plus powers ideal attaining the Hilbert function $H$,
then among all ideals with Hilbert function $H$ that also contain a
regular sequence of elements of degrees $a_1, a_2, \ldots ,a_n$, the
$\mathbb{A}$-lex plus powers ideal has the largest first graded Betti
numbers.

In light of both Bigatti and Hulett's result and Eisenbud, Green and
Harris's conjecture, the following very natural conjecture was made by
Charalambous and Evans: as long as there is an $\mathbb{A}$-lex plus
powers ideal attaining the Hilbert function $H$, then among all ideals
with Hilbert function $H$ that also contain a regular sequence of
elements of degrees $a_1, a_2, \ldots ,a_n$, the $\mathbb{A}$-lex plus
powers ideal has the largest graded Betti numbers (not just the
largest first graded Betti numbers).

As a result of Bigatti and Hulett's results, there has been much
interest in studying lex ideals.  One direction of study has led to
the introduction of $n$-type vectors by Geramita, Harima and Shin.
These $n$-type vectors are in 1-1 correspondence with Artinian lex
ideals.  Since all lex plus powers ideals are by definition Artinian,
it makes sense to look for an analogue to $n$-type vectors for lex
plus powers ideals.  We do this in section ~\ref{analogue}.  This
enables us to prove our main result quite easily: that the residual of
an $\mathbb{A}$-lex plus powers ideal in $\langle x_1^{a_1}, \ldots
,x_n^{a_n}\rangle$ is again a lex plus powers ideal.  As a consequence
of this, we show in section ~\ref{apply} that the statement that lex
plus powers ideals have largest first graded Betti numbers is
equivalent to the statement that lex plus powers ideals have largest
last graded Betti numbers (previously, it was shown in \cite{Richert}
that lex plus powers ideals having largest first graded Betti numbers
implies having the largest last graded Betti numbers; we show the
converse).

\section{Background}

Let $R=\poly{n}$ be the polynomial ring in $n$ variables over a field
$k$ with maximal ideal $m=(x_1, \dots , x_n)$, and fix an order on the
monomials, $x_1> \cdots >x_n$.  The following definition gives a
notation for referring to the degrees of the elements of a regular
sequence.
\begin{definition}
  Let $\a$ be a set of integers such that $1\le a_1 \le \cdots \le
  a_n$.  Then we call $\{f_1, \dots , f_n\}$ an $\{a_1, \dots ,
  a_n\}$-regular sequence if $\{f_1, \dots , f_n\}$ is a regular
  sequence such that $\deg(f_i)=a_i$ for $i=1, \dots ,n$.
\end{definition}

Recall that the Hilbert function $H(R/I)$ of an ideal $I$ is the
sequence $\{ \dim_k(R/I)_d\}_{d\geq 0}$.  We denote $\dim_k(R/I)_d$ by
$H(R/I,d)$.  Then given a Hilbert function $\H$, and a list of degrees
$\a$, we can compare homogeneous ideals attaining $\H$ and containing
an $\a$-regular sequence. In this comparison, we will use as a fixed
point a special ideal called an $\a$-lex plus powers ideal.

\begin{definition}[Charalambous and Evans]
  Suppose that $\A =\{a_1, \dots ,a_n\}$ is a non-decreasing list of
  integers, $a_1\ge 1$. Then a monomial ideal $L$ is a {\it lex plus
    powers ideal with respect to $\A$}, also called an {\it $\A$-lex
    plus powers ideal}, if $L$ is minimally generated by monomials
  $x_1^{a_1},\dots ,x_n^{a_n}$, $m_1, \dots ,m_l$ such that for each
  $j=1,\dots , l$, all monomials of degree $\deg(m_j)$ which are
  larger than $m_j$ in lex order are contained in $L$. We will
  abbreviate the terminology ``lex plus powers with respect to $\A$''
  by saying that $L$ is $\LA$.
\end{definition}

It is not difficult to construct (degenerative) examples of a Hilbert
function $\H$ and a list of degrees $\A=\a$ for which no $\A$-lex plus
powers ideal $L$ exists with $H(R/L)=\H$ (see \cite{Richert}). Thus we
require the following technical definition.

\begin{definition} 
  Suppose that $\H$ is a Hilbert function and $\A =\{a_1, \dots
  ,a_n\}$ is a non-decreasing list of integers, $a_1\ge 1$. We call
  $\H$ an {\it $\A$-lpp valid Hilbert function} if there exists an
  $\LA$ ideal $L$ such that $H(R/L)=\H$.  Note that if an $\LA$ ideal
  $L$ attaining a given Hilbert function $\H$ exists, then it is
  clearly unique. We will sometimes refer to this ideal as $\lppb$.
\end{definition}

Lex plus powers ideals are important because they are conjectured by
Charalambous and Evans \cite{Evans} to have extremal properties. In
order to understand in what sense lex plus powers ideals should be
extremal, we need to introduce some terminology.  Recall that the
$i,j^{\mbox{\tiny{th}}}$ graded Betti number of $I$ is defined to be
$$\beta_{i,j}^I :=(\Tor_i(R/I,k))_j.$$
We will refer to the set of all
graded Betti numbers of an ideal $I$ as $\beta^I$. It is also
convenient to make use of the notation of the computer algebra system
Macaulay 2 \cite{Macaulay2}, so we often refer to $\beta^I$ as the
Betti diagram of $I$ (the Betti diagram of $I$ is a table listing the
graded Betti numbers of $I$---counting from zero, the entry in the
$i,j^{\mbox{\tiny{th}}}$ position in this table is $\beta_{i,i+j}^I$).

\begin{definition}
  Write $\LP$ to be the set of all sets of graded Betti numbers of
  ideals $I\subset R$ containing an $\A$-regular sequence and
  attaining $\H$. Equivalently, this is the set of all Betti diagrams
  of such ideals.
\end{definition}

There is an obvious partial order on $\LP$: for
$\beta^I,\beta^J\in\LP$, we say that $\beta^I\ge \beta^J$ if
$\beta^I_{i,j}\ge \beta^J_{i,j}$ for all $i,j$. With this we can
describe the conjectured extremality of lex plus powers ideals.

\begin{conjecture}[The Lex Plus Powers conjecture] \label{c:lpp} 
  If $\H$ is $\A$-lpp valid, then writing $\lppb$ to be the $\A$-lex plus
  powers ideal attaining $\H$, $\beta^{\lppb}$ is the unique largest
  element in $\LP$.
\end{conjecture} 

There is a (on the face of it) weaker version of this conjecture due
to Eisenbud, Green, and Harris, which claims that lex plus powers
ideals should be capable of largest Hilbert function growth.
\begin{conjecture}[The Lex Plus Powers Conjecture for Hilbert
  Functions]\label{c:lpph1} Let $I\subset R$ contain an $\A$-regular
  sequence and suppose there exists an $\LA$ ideal $L$ such that
  $H(R/I,d)=H(R/L,d)$. Then $$H(R/\langle L_d\rangle,d+1)\ge
  H(R/I,d+1),$$
  where $\langle L_d\rangle$ is the ideal generated by
  the pure powers $x_1^{a_1}, \dots ,x_n^{a_n}$ and the forms in $L$
  of degree $d$.
\end{conjecture}

That the lex plus powers conjecture (LPP) implies the lex plus powers
conjecture for Hilbert functions (LPPH) is made clear by an equivalent
formulation of the latter found in \cite{Richert}:
\begin{conjecture}\label{c:lpph2}
  Given an $\A$-lpp valid Hilbert function $\H$, then
  $\beta^{\lppb}_{1,i}\ge \beta^I_{1,i}$ for all $i$ whenever
  $I\subset R$ attains $\H$ and contains an $\A$-regular sequence.
\end{conjecture} 
It is an open question whether LPPH implies LPP. Some progress was
made on this question in \cite{Richert} with the following theorem:
\begin{theorem}\label{t:socleold}
  Let $L$ be $\LA$ for some $\A=\a$ and $I$ be an ideal containing an
  $\A$-regular sequence such that $H(R/L)=H(R/I)$. If LPPH holds, then
  $\dim_k(\soc (L)_d)\ge \dim_k(\soc (I)_d)$ for all $d$.
\end{theorem}
That is, if the $\beta^{\lppb}_{1,j}$ are uniquely largest, then so
are the $\beta^{\lppb}_{n,j}$. It was not decided in that paper
whether the converse was true.  We will show in this paper that the
converse does hold. That is, we prove that the following conjecture
and LPPH are equivalent:
\begin{conjecture}\label{c:socle}
  Let $L$ be $\LA$ for some $\A=\a$ and $I$ be an ideal containing an
  $\A$-regular sequence such that $H(R/L)=H(R/I)$. Then
  $\beta^{L}_{n,j}\ge \beta^I_{n,j}$, that is, $\dim_k(\soc (L)_d)\ge
  \dim_k(\soc (I)_d)$ for all $d$.
\end{conjecture}

This result will be a natural application of our main result, where we
show that the residual of an $LPP(\mathbb{A})$-ideal in $\langle
x_1^{b_1}, \ldots ,x_n^{b_n}\rangle$ where $a_i\leq b_i$ for all $i$
is again a lex plus powers ideal.

We recall here one further theorem, a result of Stanley.

\begin{theorem}[Stanley]\label{t:stanley}
  For every $R$-module $M$,
\begin{displaymath}
  \sum_{d=0}^\infty H(M,d)t^d=\frac{\sum_{d=0}^\infty \sum_{i=0}^{n}
    (-1)^i\beta^M_{i,d}t^d}{(1-t)^n}.
\end{displaymath}
\end{theorem} 

This theorem simply states that fixing a Hilbert function fixes the
alternating sum of the graded Betti numbers of any ideal attaining it.
In particular, if $I$ and $J$ have $H(R/I)=H(R/J)$, then
$\sum_{i=0}^n(-1)^i\beta^I_{i,j}=\sum_{i=0}^n(-1)^i\beta^J_{i,j}$ for
all $j$.  This implies that for $\rho$ the regularity of
$H(R/I)$, $\beta^I_{n,\rho+n}=\beta^J_{n,\rho+n}$ and
$\beta^I_{n-1,\rho+n-1}-\beta^I_{n,\rho+n-1}=
\beta^J_{n-1,\rho+n-1}-\beta^J_{n,\rho+n-1}$. These last two facts
will prove useful in section \ref{apply}.

\section{The Hilbert function of lex plus powers ideals}

In this section, we state a characterization of the Hilbert functions
which can occur for $\a$-lex plus powers ideals.  This
characterization follows from the work of Clements and Lindstrom and
will be useful in the next section when we find an alternative to the
Hilbert functions of lex plus powers ideals similar to the $n$-type
vectors found by Geramita, Harima and Shin in \cite{ghs} for Hilbert
functions of lex ideals.  For more details than provided here on the
relationship between the work of Clements and Lindstrom and Macaulay's
O-sequences, see \cite{cr}.

\begin{definition}
Let $\A=\a$. Then a lex plus powers ideal $L$ is said to be {\em lex
    plus powers with respect to $\le \A$}, or $\LL{\A}$, if $L$
  contains the $\A$-regular sequence $\{x_1^{a_1}, \dots ,
  x_n^{a_n}\}$. Note that a $\B=\b$-lex plus powers ideal is $\LL{\A}$
  if and only if $\B\le \A$, that is, if $b_i\le a_i$ for all $i=1,
  \dots , n$.
\end{definition}

Although Clements and Lindstrom used different terminology, the following special case of the LPPH conjecture can be found in their paper \cite{cl}.

\begin{theorem} Let $\A=\a$, $L$ be $LPP(\A )$ and $I$ be any monomial ideal in $R=k[x_1,\ldots ,x_n]$ containing $x_1^{a_1},\ldots ,x_n^{a_n}$ such that $H(R/I,d)=H(R/L,d)$.  Then $H(R/I,d+1)\leq H(R/\langle L_d\rangle ,d+1)$.
\end{theorem}

Since any $\LL{\A}$-ideal is a monomial ideal containing $x_1^{a_1},\ldots ,x_n^{a_n}$, we obtain:

\begin{corollary} \label{candl} Let $\A=\a$ and $L$ be $LPP(\A )$ and $I$ be an $\LL{\A}$ ideal such that $H(R/L,d)=H(R/I,d)$. Then $H(R/I,d+1)\leq H(R/\langle L_d\rangle ,d+1)$.
\end{corollary}

Keeping in the Macaulayesque mindset, we introduce the following notation:

\begin{definition} \label{mac} 
  Let $\A=\a$.  Let $L$ be an $LPP(\A )$-ideal satisfying
  $H(R/L,d)=h$.  Then define $h^{{\langle d\rangle}_{\A}}:=
  H(R/\langle L_d\rangle ,d+1)$.  Furthermore, let $S=\{ c_i\}_{i\geq
    0}$ be a sequence satisfying $c_0=1$ and $c_{i+1}\leq
  c_i^{{\langle i\rangle}_{\A}}$ for all $i$.  Then $S$ is said to be
  an $\LL{\A}$-sequence.
\end{definition}

\begin{remark} 
  In the notation of Definition ~\ref{mac}, Corollary ~\ref{candl}
  says that $H$ is the Hilbert function of an $\LL{\A}$-ideal if and
  only if $H$ is an $\LL{\A}$-sequence.  Note that, to determine the
  Hilbert functions of $LPP(\A )$-ideals, we cannot simply eliminate
  the sequences that are $lpp_{\leq}(\B )$-sequences for $\B \leq\A$,
  but $\B\neq \A$ from the set of $lpp_{\leq}{(\A )}$-sequences.  This
  is because of the possibility of overlap.  For example, $I =\langle
  x^2,y^3,z^4,xy^2,xyz,xz^2,y^2z^2\rangle$ and $J=\langle
  x^2,y^3,z^3,xy^2,xyz\rangle$ are respectively $LPP(\{ 2,3,4\} )$ and
  $LPP(\{ 2,3,3\} )$-ideals, both having Hilbert function $H$ = 1 3 5
  1 0 $\rightarrow$.
\end{remark}

Greene and Kleitman (\cite{comb}) found a Macaulayesque way of
describing $h^{\langle i\rangle _{\mathbb{A}}}$, which we wish to
consider in some detail, since we will be using their notation in
later parts of this paper.  Before doing so, we recall Macaulay's
methods.

  Let $d,h\in\N$ be given. Then it is well known that there are unique
  integers $k(d)>k(d-1)> \cdots > k(1)\ge 0$ such that $h={k(d)\choose
    d}+{k(d-1)\choose d-1} + \cdots + {k(1)\choose 1}$. Macaulay's
  theorem states that if $h$ is the value of the Hilbert function of a
  graded module in degree $d$, then $H(M,d+1)\le {k(d)+1\choose
    d+1}+{k(d-1)+1\choose d-1+1} + \cdots + {k(1)+1\choose 1+1}$, and
  this bound is sharp.  The process of obtaining the $k(i)$ and computing the bound can be
  beautifully visualized by writing Pascal's triangle as a rectangle:

$$
\begin{array}{rrrrrrr}     
                           0&1&2&3&4&5&  \\  
                            & & & & & &  \\
                           1&1&1&1&1&1&\ldots\\
                           1&2&3&4&5&6&\\
                           1&3&6&10&15&21&\\
                           1&4&10&20&35&56&\\
                           1&5&15&35&70&126&\\
                           1&6&21&56&126&252&\\
                           \vdots& & & & & &\ddots\\
                            & & & & & &  
\end{array}
$$
 
\begin{example}
Suppose that $M$ is a graded module such that $H(M,3)=32$. Then to
obtain an upper bound for $H(M,4)$, one must first find the $k(i)$
which uniquely describe $32$ in degree $3$. First, look at the column
numbered 3, and pick the largest number that is at most 32, namely 20.
This is $3$ rows down from the top, so we take $k(3)=3+3=6$. Then
look at the column numbered 2 and pick the largest number that is at most 
$32-20=12$, namely 10. This is again $3$ rows down, so we take
$k(2)=2+3=5$. Finally, pick the 2 from the column numbered 1, which is
$1$ row down, so we take $k(1)=1+1=2$. 
$$
\begin{array}{rrrrrrr}     
                           0&1&2&3&4&5&  \\  
                            & & & & & &  \\
                           1&1&1&1&1&1&\ldots\\
                           1&\frame{2}& 3&4&5&6&\\
                           1&3&6&10&15&21&\\
                           1&4&\frame{10}& \frame{20}& 35&56&\\
                           1&5&15&35&70&126&\\
                           1&6&21&56&126&252&\\
                           \vdots& & & & & &\ddots\\
                            & & & & & &  
\end{array}
$$
Recalling that the number in the $i$th row and $j$th column of
Pascal's triangle is ${i+j-1\choose j}$, it is evident that we have
just found $32={6\choose 3}+{5\choose 2}+{2\choose 1}$ (note that
${6\choose 3}=20$, ${5\choose 2}=10$, and ${2\choose 1}=2$). Then to
compute the bound for $H(M,4)$, we need ${6+1\choose 3+1}+{5+1\choose
  2+1}+{2+1\choose 1+1}$, and this is obtained by taking the number
one column to the right of each of the boxed integers in the rectangular version of Pascal's
triangle:
$$
\begin{array}{rrrrrrr}     
                           0&1&2&3&4&5&  \\  
                            & & & & & &  \\
                           1&1&1&1&1&1&\ldots\\
                           1&\frame{2}&\rightarrow 3&4&5&6&\\
                           1&3&6&10&15&21&\\
                           1&4&\frame{10}&\rightarrow
                            \frame{20}&\rightarrow 35&56&\\ 
                           1&5&15&35&70&126&\\
                           1&6&21&56&126&252&\\
                           \vdots& & & & & &\ddots\\
                            & & & & & &  
\end{array}
$$
The result is $H(M,4)\le 35+20+3=58$.

\end{example}

\begin{remark}
  There is a precise relationship between monomials of degree $i$ and
  $i$-binomial expansions.  Namely, if $h=\binom{m_i}{i} + \binom{m_{i
      -1}}{i-1} +\ldots +\binom{m_j}{j}$, then $h$ is the codimension
  of a lex-segment in degree $i$ in the polynomial ring in $n=m_i-i+2$
  variables.  Letting $m$ be the smallest monomial of degree $i$ in
  this lex-segment, we associate $h$ to $m$.  Namely, let $\alpha _r
  =\#\{ t|m_t-t=n-1-r\}$ for $1\leq r\leq n-1$.  Then the lex segment
  ending in the monomial $m=x_1^{\alpha _1}x_2^{\alpha _2}\ldots
  x_{n-1}^{\alpha _{n-1}}x_n^{i-(\alpha _1+\cdots +\alpha _{n-1})}$
  has codimension $h$ in $k[x_1,\ldots ,x_n]$. (See \cite{rob} for
  details.)

Since $\alpha _1 +\alpha _2 +\ldots +\alpha _{n-1}$ is the number of
terms in the $i$-binomial expansion of $h$, we see that $i-(\alpha _1
+\alpha _2 +\ldots +\alpha _{n-1} )=j-1$, so we can rewrite $m$ as
$x_1^{\alpha _1}x_2^{\alpha _2}\ldots
x_{n-1}^{\alpha_{n-1}}x_n^{j-1}$.  In fact, this correspondence could
have been used to define $i$-binomial expansions in the first place,
and is the reason why they are so valuable in the study of Hilbert
functions.

\end{remark}

We now wish to state the growth bound for $\LL{\A}$ ideals in terms of the notation used by Greene and Kleitman in \cite{comb}.  Let $d_1\leq d_2\leq \ldots \leq d_k$ and put $e_1 :=d_k -1, e_2:= d_{k-1} -1, \ldots ,e_k := d_1 -1$.  Then they used the notation $\binom{e_1,\ldots ,e_k}{i}$ to be $\Delta H(R/I,i)$, where $\Delta$ represents the first difference function and $I\subset k[x_0,\ldots ,x_n]$ is the ideal of a complete intersection of type $(d_1, \ldots ,d_k)$. Note that $\binom{e_1}{i}$ is not the usual binomial coefficient; $\binom{e_1}{i}$ is 1 if $0\leq i\leq e_1$ and is 0 if $i>e_1$. This will allow us to state the LPPH conjecture using their Macaulayesque form, but first, we need a result stated in \cite{comb}.  

 \begin{defprop} Let $\mathbb{A}=\a$ and $d$ be given and let $0< h\leq H(R/(x_1^{a_1},\cdots ,x_n^{a_n}),d)$.  Let $a'_i=a_i-1$.  Then $h$ can be written uniquely in the form 
   
   {\tiny $$h=\binom{a'_n,a'_{n-1},\cdots
       ,a'_{n-(k(d)-d)}}{d}+\binom{a'_n,a'_{n-1},\cdots
       ,a'_{n-(k(d-1)-(d-1))}}{d-1}+\cdots
     +\binom{a'_n,a'_{n-1},\cdots ,a'_{n-(k(j)-j)}}{j}.$$} where
   $k(d)>k(d-1)>\ldots >k(j)\geq j\geq 1$ and $\# \{ t|k(t)-t = i\} <
   a_ {n-i-1}$ and the last term is non-zero.
   \\
   \\
   We refer to this expression as the {\bf $d_{\mathbb{A}}$-Macaulay
     expansion} for $k$. Furthermore, {\tiny $$h^{\langle d\rangle
       _{\A}} := \binom{a'_n,a'_{n-1},\cdots
       ,a'_{n-(k(d)-d)}}{d+1}+\binom{a'_n,a'_{n-1},\cdots
       ,a'_{n-(k(d-1)-(d-1))}}{d}+\cdots +\binom{a'_n,a'_{n-1},\cdots
       ,a'_{n-(k(j)-j)}}{j+1}.$$}
\end{defprop}

One way to look at this proposition is through the correspondence between monomials $m$ and the codimension of the lex-segments ending in monomial $m$. Given a monomial $m=x_1^{\alpha_1}\cdots x_n^{\alpha _n}$, write the expansion for which $\alpha _i =\#\{ t|k(t)-t = n-1-i\}$ for $1\leq i\leq n-1$ and then remove any zero terms at the end.

\begin{example} Let $R=k[x_1,x_2,x_3]$ and $\mathbb{A} = \{ 3, 4, 11\}$. The monomials of degree 12 in $k[x_1,x_2,x_3]/\langle x_1^3,x_2^4,x_3^{11}\rangle$ with their codimensions are listed below:
$$
\begin{array}{ll}
x_1^2x_2^3x_3^7 & \binom{10,3}{12}+\binom{10,3}{11}+\binom{10}{10}+\binom{10}{9}+\binom{10}{8}=8\\
x_1^2x_2^2x_3^8 & \binom{10,3}{12}+\binom{10,3}{11}+\binom{10}{10}+\binom{10}{9}=7\\
x_1^2x_2x_3^9 & \binom{10,3}{12}+\binom{10,3}{11}+\binom{10}{10}=6\\
x_1^2x_3^{10} & \binom{10,3}{12}+\binom{10,3}{11}=5\\
x_1x_2^3x_3^8 & \binom{10,3}{12}+\binom{10}{11}+\binom{10}{10}+\binom{10}{9}=4\\
x_1x_2^2x_3^9 & \binom{10,3}{12}+\binom{10}{11}+\binom{10}{10}=3\\
x_1x_2x_3^{10} & \binom{10,3}{12}+\binom{10}{11}=\binom{10,3}{12}=2\\
x_2^3x_3^{9} & \binom{10}{12}+\binom{10}{11}+\binom{10}{10}=1\\
x_2^2x_3^{10} & \binom{10}{12}+\binom{10}{11}=0
\end{array}
$$
\end{example}

\begin{conjecture} (Restatement of the LPPH Conjecture): Let $I\subset R$ contain an ${\A}$-regular sequence and suppose there exists an $LPP(\mathbb{A})$-ideal $L$ such that $H(R/I,d) = H(R/L,d)$.  Then $H(R/I,d+1)\leq H(R/I,d)^{\langle d\rangle_{\A}}$.
\end{conjecture}

\begin{example}
Suppose for instance that $R=k[x_1,x_2,x_3]$, $\A=\{3,4,11\}$, and $L$
is $\LL{A}$ with $H(R/L,4)=10$. Then we consider the following
rectangle
\begin{tiny}
$$
\begin{array}{rrrrrrrrrrrrrrrrrrr}     
         &0&1&2&3&4 &5 &6 &7 &8 &9 &10&11&12&13&14&15&16& \ldots\\  
         & & & & &  &  &  &  &  &  &  &  &  &  &  &  &  &   \\
(1,1,11):&1&1&1&1&1 &1 &1 &1 &1 &1 &1 &0 &\rightarrow& & & & & \\
(1,4,11):& 1&2&3&4&4 &4 &4 &4 &4 &4 &4 &3 & 2& 1&0 &\rightarrow& & \\
(3,4,11):& 1&3&6&9&11&12&12&12&12&12&12&11&9 &6 &3 &1 &0 &\rightarrow 
\end{array}
$$
\end{tiny}
where we have written $(a_1,a_2,a_3)$ beside the row that consists of $\Delta H(R/I)$ for $I$ a complete intersection of type $(a_1,a_2,a_3)$.  The top row is thus $\binom{10}{i}$ for $i\geq 0$, the second row is $\binom{10,3}{i}$ for $i\geq 0$, and the third row is $\binom{10,3,2}{i}$ for $i\geq 0$.

The largest number in the column numbered 4
which is at most 10 is 4. In the column numbered 3, we take the largest
number that is at most $10-4=6$, which is 4. In the column numbered 2, we
take 1, and finally in the column numbered 1, we pick  1. This expresses 
$10$ as a $4_{\A}$-Macaulay expansion:
\begin{tiny}
$$
\begin{array}{rrrrrrrrrrrrrrrrrrr}     
         &0&1&2&3&4 &5 &6 &7 &8 &9 &10&11&12&13&14&15&16& \ldots\\  
         & & & & &  &  &  &  &  &  &  &  &  &  &  &  &  &   \\
(1,1,11):&1&\frame{1}& \frame{1}& 1&1 &1 &1 &1
         &1 &1 &1 &0 && & & & & \\ 
(1,4,11):& 1&2&3&\frame{4}& \frame{4}& 4 &4 &4
         &4 &4 &4 &3 & 2& 1&0 && & \\ 
(3,4,11):& 1&3&6&9&11&12&12&12&12&12&12&11&9 &6 &3 &1 &0 &
\end{array}
$$
\end{tiny}
Note that the number to the right of ${e_1,\ldots ,e_k\choose i}$ is just ${e_1,\ldots ,e_k\choose i+1}$.  Thus, to calculate $10^{\langle 4\rangle _{\A}}$, the bound for $H(R/L,5)$, we again sum the numbers to the right of our boxed integers.  
\begin{tiny}
$$
\begin{array}{rrrrrrrrrrrrrrrrrrr}     
         &0&1&2&3&4 &5 &6 &7 &8 &9 &10&11&12&13&14&15&16& \ldots\\  
         & & & & &  &  &  &  &  &  &  &  &  &  &  &  &  &   \\
(1,1,11):&1&\frame{1}&\rightarrow \frame{1}&\rightarrow 1&1 &1 &1 &1
         &1 &1 &1 &0 &\rightarrow& & & & & \\ 
(1,4,11):& 1&2&3&\frame{4}&\rightarrow \frame{4}&\rightarrow 4 &4 &4
         &4 &4 &4 &3 & 2& 1&0 &\rightarrow& & \\ 
(3,4,11):& 1&3&6&9&11&12&12&12&12&12&12&11&9 &6 &3 &1 &0 &\rightarrow 
\end{array}
$$
\end{tiny}
Thus, we find that $H(R/L,5)\le 4+4+1+1=10$.
\end{example}

\begin{example} 
  Suppose that $L$ is an $\A=\{3,4,11\}$ lex plus powers ideal and
  $H(R/L,12)=7$. The monomials of degree 12 not in $I$ are $$x^2yz^9,\ 
  x^2z^{10},\ xy^3z^8,\ xy^2z^9,\ xyz^{10},\ y^3z^9,\ \mbox{and}\ 
  y^2z^{10}$$
  and so in degree 13, at most the following monomials are
  not in $I$:
  $$x^2yz^{10},\ xy^3z^9,\ xy^2z^{10},\ \mbox{and}\ y^3z^{10}.$$
  Then
  the diagram looks like
\begin{tiny}
$$
\begin{array}{rrrrrrrrrrrrrrrrrrr}     
         &0&1&2&3&4 &5 &6 &7 &8 &9 &10&11&12&13&14&15&16& \ldots\\  
         & & & & &  &  &  &  &  &  &  &  &  &  &  &  &  &   \\
(1,1,11):&1&1&1&1&1 &1 &1 &1 &1 &\frame{1} &\rightarrow \frame{1}
         &\rightarrow 0 &\rightarrow& & & & & \\ 
(1,4,11):& 1&2&3&4&4 &4 &4 &4 &4 &4 &4 &\frame{3} &\rightarrow
         \frame{2}&\rightarrow 1&0 &\rightarrow& & \\ 
(3,4,11):& 1&3&6&9&11&12&12&12&12&12&12&11&9 &6 &3 &1 &0 &\rightarrow 
\end{array}
$$
\end{tiny}
so that as expected $H(R/L,13)\le 1+2+0+1=4$. 

\end{example}

\section{An analogue to $n$-type vectors for lex plus powers ideals.}\label{analogue}
We wish to define a vector that will correspond in a natural way to lex plus powers ideals.  This will be an analogue to the $n$-type vectors that correspond to lex ideals. Let $a\leq b$. Then any $LPP(a,b)$-ideal is of the form
$$L=\langle x^a, x^{a-1}y^{d_1}, x^{a-2}y^{d_2}, \ldots ,x^{a-s}y^{d_s},y^b\rangle$$
where $d_1<d_2<\ldots <d_s<b$. We associate to $L$ the vector $\T =(d_1, d_2,\ldots ,d_s,b,\ldots ,b)$ where there are $a-s$ $b$'s and $a\leq b$.  The condition that $a\leq b$ is crucial, for otherwise the ideal would not be lex plus powers.

\begin{example} If we put $\T = (2,4,5,5,5,5)$, the associated ideal would be $I=\langle x^6, x^5y^2, x^4y^4,y^5\rangle$.  Since this violates the condition that the powers of the variables be in non-decreasing order, the ideal is not $LPP(5,6)$.  The $LPP(5,6)$-ideal with the same Hilbert function as $I$ is $J=\langle x^5,x^4y^3,x^3y^5,y^6\rangle$ and this corresponds to the vector (3,5,6,6,6).  They both have the same graded Betti numbers, but for uniqueness purposes, we choose $J$ as the $LPP(5,6)$-ideal. 
\end{example} 

\begin{remark}
In three variables, it is easy to construct (\cite[Remark 4.3]{res}) many ideals which satisfy all the requirements of lex plus powers ideals except the condition that the powers of the variables are in non-decreasing order, and do not actually have the same graded Betti numbers as the lex plus powers ideal.  
\end{remark}

\begin{definition} Let $\A =\{ a_1,\ldots ,a_n\}$.

If $n=1$ and $\T = (d)$ for some $d\leq a_1$, we say that $\T$ is an $\LL{\A}$-vector.  We say that $\T = \T_{c.i.(\A )}$ if $\T =(a_1)$.  We put
$\sigma (\T )=l(\T )=\alpha _{\mathbb{A}}(\T )=d$ unless $\T =\T_{c.i.(\A )}$, in which
case we put $l(\T )=\sigma (\T )=a_1$ and $\alpha _{\mathbb{A}}(\T )=\infty$.

If $n>1$, then $\T =(\T _1, \ldots ,\T _u)$ is an $\LL{\bb{A}}$-vector if the following conditions all hold:
$u\leq a_1$, $u \leq l(\T _u)$, each $\T _i$ is an $\LL{\bb{A}_2}$-vector (in particular, $l(\T _u)\leq a_2$) and $\sigma (\T _i)<\alpha _{\A _2}(\T _{i+1})$ for $1\leq i \leq u-1$.

We define $l(\T )=u$ to be the length of $\T$, and $\sigma (\T )$ and $\alpha _{\mathbb{A}}(\T )$ as follows:
$$
\begin{array}{lll}
\sigma (\T )&=&\left\{ \begin{array}{ll}
                      \sigma (\T _u)&\mbox{if }\T _u\neq \T_{c.i.(\A _2)}\\
                      \sigma (\T _u) +s-1&\mbox{if }\T _u = \T_{c.i.(\A _2)}\mbox{ where } s=\# i\ s.t.\ \T _i =\T _u.
                      \end{array}\right. \\
\\
\alpha _{\mathbb{A}}(\T )&=&\left\{ \begin{array}{ll}
                      l(\T )&\mbox{if }l(\T )<a_1\\
                      l(\T ) +\alpha _{\mathbb{A}_2}(\T _1)-1&\mbox{if }l(\T ) = a_1.
                      \end{array}\right. 
\end{array}
$$     
Finally, we say that $\T =\T_{c.i.(\A )}$ if $l(\T )=a_1$
  and $\T _i =\T _{c.i.(\A _2)}$ for each $i$.
\end{definition}

\begin{remark} $\alpha _{\mathbb{A}}(\T )<\infty$ unless $\T =\T_{c.i.(\A )}$.  Furthermore, $\alpha _{\mathbb{A}}(\T)\leq \sigma (\T )$ unless $\T =\T _{c.i.(\A )}$.
\end{remark}

{\bf Notation:} For convenience, we will denote the vector $((d_1),\ldots ,(d_m))$ by $(d_1,\ldots ,d_m)$.  Thus, for example, the vector ((1),(3),(4)) will be written as (1,3,4), and the vector (((1),(2)),((1),(3),(4))) will be written as ((1,2),(1,3,4)).  This does however create confusion since $(d_1)$ could denote either the vector $((d_1))$ or the vector $(d_1)$.  If there is ever any confusion, we will explicitly state what we are referring to.

\begin{example} Let 
$$\T =(\T _1, \T _2,\T _3, \T _4, \T _5) = ((1,2),(1,3,4),(2,3,6,6),(5,6,6,6),(6,6,6,6)),$$
where each $\T _i$ is an $lpp_{\leq}({4,6})$-vector.  Then both $(\T _1, \T _2,\T _3, \T _4)$ and $(\T _2,\T _3, \T _4, \T _5)$ are $\LL{\mathbb{A}}$-vectors where $\mathbb{A}= \mathbb\{ 4,4,6\} $ since $\sigma (\T _1)=2<\alpha _{\bb{A}_2}(\T _2)=3$, $\sigma (\T _2)=4<\alpha _{\bb{A}_2}(\T _3)=4+2-1=5$, $\sigma (\T _3)=6+2-1=7<\alpha _{\bb{A}_2}(\T _4)=4+5-1=8$ and $\sigma (\T _4)=6+3-1=8<\alpha _{\bb{A}_2}(\T _5)=\infty$.  However, $\T$ is not an $\LL{\mathbb{A}}$-vector for any $\mathbb{A}=\{ a_1,a_2,a_3\}$, for suppose it were.  Then $a_1\geq l(\T )=5$.  Since $\mathbb{A} =\{ a_1,a_2,a_3\}$ must satisfy $a_1\leq a_2\leq a_3$, we also have $a_2\geq 5$.  Then $\alpha _{\mathbb{A}_2}(\T _5)=4$ and $\sigma (\T _4)=8$, contradicting that $\sigma (\T_4)<\alpha_{\mathbb{A}_2} (\T_5).$ Notice also that $\T_3, \T_4$ and $\T_5$ are all $lpp_{\leq}({5,6})$-vectors, but are not $lpp_{\leq}({4,7})$-vectors.
\end{example}

To an $\LL{\bb{A}}$-vector $\T$, it is natural to associate an ideal $W_{\T}$ as follows:

\begin{definition}
If $n=1$ (so that $\A =\{ a_1\}$) and $\T$ is an $\LL{\bb{A}}$-vector, say $\T=(d)$ with $d\leq a_1$, then define $W_{\T} := \langle x_1^d\rangle$ in $k[x_1]$.

If $n>1$ and $\T$ is an $\LL{\bb{A}}$-vector, say $\T =(\T _1,\ldots ,\T _u)$ with $u\leq a_1$, then define 
$$W_{\T} := \langle x_1^u, x_1^{u-1}\overline{W_{\T _1}},\ldots ,x_1\overline{W_{\T _{u-1}}},\overline{W_{\T _u}}\rangle$$
where $\overline{W_{\T _i}}$ is the image in $k[x_2,\ldots ,x_n]$ under the isomorphism induced by $x_i\to x_{i+1}$ of the ideal $W_{\T _i}\subset k[x_1,\ldots ,x_{n-1}]$ obtained by induction.   
\end{definition}

\begin{remark} If $\T =\T _{c.i.(\A )}$, then $W_{\T} =\langle x_1^{a_1}, \ldots ,x_n^{a_n}\rangle$.  To see this, note that if $n=1$ and $\T =(a_1)$, then $W_{\T} =\langle x_1^{a_1}\rangle$ and by induction, if $\T =(\T _{c.i.(\A _2)}, \ldots ,\T _{c.i.(\A _2)})$ with $l(\T )=a_1$, then $W_{\T} =\langle x_1^{a_1}, \overline{W_{\T _{c.i.(\A _2)}}}\rangle = \langle x_1^{a_1}, \ldots ,x_n^{a_n}\rangle$, as required.
\end{remark}

Before showing that $W_{\T}$ is an $\LL{\bb{A}}$-ideal if
$\T$ is an $\LL{\bb{A}}$-vector, we first show that $\alpha (\T)$ is the
smallest degree of any element of $W_{\T}$ not in $\langle x_1^{a_1},
\ldots ,x_n^{a_n}\rangle$ and that $\sigma (\T )-1$ is the largest degree of any element of $k[x_1,\ldots ,x_n]$ not in $W_{\T}$.  In fact, we give names to these parameters for any ideal containing $\langle x_1^{a_1}, \ldots ,x_n^{a_n}\rangle$.  

\begin{definition} Let $I$ be any ideal of $k[x_1,\ldots ,x_n]$ containing $\langle x_1^{a_1}, \ldots ,x_n^{a_n}\rangle$.  Then put 
$$\alpha_{\mathbb{A}}(I)=\min \{ i|f\in I\setminus \langle x_1^{a_1}, \ldots ,x_n^{a_n}\rangle, \deg f=i\} \hspace{.5cm}\mbox{and}$$
$$\sigma (I)=\min \{ i|I_i =k[x_1,\ldots ,x_n]_i\} .$$
\end{definition} 

We use $\alpha _{\mathbb{A}}$ instead of $\alpha$ to distinguish it from the usual $\alpha$, which is just $\alpha (I)=\min \{ i|f\in I, \deg f=i\}$.  $\sigma (I)$ is defined as usual.

\begin{lemma} \label{alphalem} Let $\T $ be an $\LL{\bb{A}}$-vector.  Then $\alpha _{\bb{A}}(W_{\T} )=\alpha _{\bb{A}}(\T )$.  
\end{lemma}

\begin{proof} The result is clear for $n=1$, so assume that $n>1$.
  Furthermore, the result is clear if $\T =\T _{c.i.(\A )}$, so we assume this is not the case.
  
  Let $\T = (\T _1, \ldots ,\T _u, \T _u, \ldots ,\T _u)$, where $l(\T
  ) =u+v$, so there are $v+1$ $\T _u$'s.  Then $W_{\T} = \langle
  x_1^{u+v}, x_1^{u+v -1}\overline{W_{\T _1}}, \ldots
  ,x_1^{v+1}\overline{W_{\T _{u-1}}}, \overline{W_{\T _u}}\rangle$.
  There are four cases to consider, determined by whether or not
  $\T_u=\T_{c.i.(A)}$ and whether or not $u+v=a_1$. Each proof is
  similar, so we include only the case for which $\T_u=\T_{c.i.(A)}$
  and $u+v=a_1$ as a representative.

We know by the induction hypothesis that the smallest degree of any
element of $\overline{W_{\T _i}}$ not in $\langle x_2^{a_2}, \ldots
,x_n^{a_n}\rangle$ is $\alpha _{\bb{A}_2}(\T _i)$.  Now,
$\overline{W_{\T _u}}=\langle x_2^{a_2}, \ldots ,x_n^{a_n} \rangle$,
so we can ignore it.  Now for $i<u$, we have $\alpha _{\bb{A}_2}(\T
_i)\leq \sigma (\T _i)<\alpha _{\bb{A}_2}(\T _{i+1})$, so

$\alpha _{\A }(W_\T  )=u+v-1+\alpha _{\bb{A}_2}(\T _1) =a_1 -1 +\alpha _{\bb{A}_2}(\T _1)=l(\T ) +\alpha _{\bb{A}_2}(\T _1) -1=\alpha _{\bb{A}}(\T ).$

\end{proof}

\begin{lemma} \label{sigmalem} Let $\T $ be an $\LL{\bb{A}}$-vector.  Then $\sigma (W_{\T} )=\sigma (\T )$.
\end{lemma}

\begin{proof} If $n=1$, the result is clear, so suppose that $n>1$.  Let $\T =(\T _1, \ldots ,\T _u, \ldots ,\T _u)$, where $l(\T )=u+v$ and there are $v+1$ $\T _u$'s (if $v>0$ then $\T_u$ is necessarily $\T_{c.i.(\mathbb{A}_2)}$).  Then we have $W_{\T} = \langle x_1^{u+v}, x_1^{u+v-1}\overline{W_{\T _1}}, \ldots ,x_1^{v+1}\overline{W_{\T _{u-1}}}, \overline{W_{\T _u}}\rangle$.  We know that there is an element of $x_1^vk[x_2, \ldots ,x_n]_{\sigma (\T _u)-1}$ that is not in $W_{\T }$.  We claim that $(W_{\T})_{\sigma (\T _u) +v}=k[x_1, \ldots ,x_n] _{\sigma (\T _u) +v}$.  So let $f$ be a monomial of degree $\sigma (\T _u) +v$.  If $x_1^{v+1}|f$, then we have that $f\in x_1^{v+i}k[x_2, \ldots ,x_n]_{\sigma (\T _u)-i}$ for some $i$.  But $\sigma (\T _u)-i\geq \sigma (\T _{u-i})$, so $f\in W_{\T }$.  If $x_1^{v+1}$ does not divide $f$, then the part of $f$ in $k[x_2, \ldots ,x_n]$ has degree at least $\sigma (\T _u)$, so $f\in W_{\T}$, as required.

\end{proof}

\begin{theorem} If $\T$ is an $\LL{\bb{A}}$-vector, then $W_{\T}$ is an $\LL{\bb{A}}$-ideal.  
\end{theorem}

\begin{proof} If $n=1$, the result is clear.  So, let $\T = (\T _1, \ldots ,\T _u, \T _u, \ldots ,\T _u),$ where $l(\T )=u+v$, so there are $v+1$ $\T _u$'s.  Then $$W_{\T} =
  \langle x_1^{u+v}, x_1^{u+v -1}\overline{W_{\T _1}}, \ldots
  ,x_1^{v+1}\overline{W_{\T _{u-1}}}, \overline{W_{\T _u}}\rangle .$$

By the induction hypothesis, each $W_{\T _i}$ is an
$\LL{\bb{A}_2}$-ideal.  Furthermore, since $l(\T)\leq a_1$, and
$l(\T)\leq l(\T _u)\leq a_2$, it is enough to show that any largest
degree element of $x_1^{u+v-i}k[x_2,\ldots ,x_n]$ not in
$x_1^{u+v-i}\overline{W_{\T _i}}$ has degree smaller than any smallest
degree element of $x_1^{u+v-(i+1)}\overline{W_{\T_{i+1}}}$ not in
$\langle x_1^{a_1},x_2^{a_2}, \ldots ,x_n^{a_n}\rangle$.  Thus, we
need to show that $\sigma (\overline{W_{\T _i}})-1+u+v-i<\alpha _{\A
  _2}(\overline{W_{\T _{i+1}}}) + u+v-(i+1)$ or in other words (from
Lemmas ~\ref{alphalem} and ~\ref{sigmalem}) that $\sigma (\T
_i)<\alpha_{\A _2}(\T _{i+1})$.  Since $\T$ is an $\LL{\A}$-vector, we
are done.
\end{proof}

To a given $\LL{\bb{A}}$-vector, we associate a Hilbert function as follows:

\begin{definition}
If $n=1$, so that $\T =(d)$ is an $\LL{\bb{A}}$-vector, then define $H_{\T}$ to be the sequence $H_{\T} := 1\ 1\ 1\ \ldots 1\ 0\rightarrow$ with $d$ 1's.

If $\T =(\T _1, \ldots ,\T _u)$, then define $H_{\T}$ to be the sequence 
$$H_{\T} (i) := \sum_{j=1}^u H_{\T _j} (i-u+j).$$
\end{definition}

We want to show that if $\T$ is an $\LL{\bb{A}}$-vector, then $H(R/W_{\T})=H_{\T}$.  We need the following lemmas.

\begin{lemma}\label{alphaless} Let $\T$ be an $\LL{\A}$-vector.  Then $\alpha _{\A}(\T )\leq \alpha _{\A _2}(\T_{l(\T )})$.
\end{lemma}

\begin{proof} The proof is easy and hence omitted. 

\end{proof}

\begin{lemma}\label{sigmaless} Let $\T$ be an $\LL{\A}$-vector.  Let $0\leq j\leq l(\T )-1$.  Then $\sigma (\T )-j\geq \sigma (\T _{l(\T )-j})$.  
\end{lemma}

\begin{proof} The proof is easy and hence omitted. 
\end{proof}
        
\begin{lemma} \label{containment}Let $\T =(\T _1,\ldots ,\T _u,\T _u, \ldots ,\T _u)$ be an $\LL{\mathbb{A}}$-type vector.  Then $W_{\T _i}\supsetneq W_{\T _{i+1}}$ for all $i=1,\ldots ,u-1$.
\end{lemma}

\begin{proof} For notational convenience, we leave out the bar notation and assume it to be understood, so we write $\overline{W_{\T_1}}$ as $W_{\T_1}$ and $\overline{(\overline{W_{\T_1}})_1}$ as $(W_{\T_1})_1$.

We use induction on $n$, where $n$ is the length of $\mathbb{A}$.   

$n=2$:  $\T =(e_1, \ldots ,e_u,e_u,\ldots ,e_u)$.  We need to show that $\langle x_1^{e_i}\rangle\supsetneq \langle x_1^{e_{i+1}}\rangle$ for $i<u$, but this is true since $e_{i+1}>e_i$.

$n>2$: We first show that $((\T _i )_{l(\T _i)-j},(\T _{i+1})_{l(\T _{i+1})-j})$ is an $lpp_{\leq}(\mathbb{A}_2)$-type vector for $0\leq j\leq l(\T _i)-1$. Let $\T _i =((\T _i )_1,(\T _i)_2,\ldots ,(\T _i)_{l(\T _i)})$ and $\T _{i+1} =((\T _{i+1} )_1,\ldots ,(\T _{i+1})_{l(\T _{i+1})}).$  Now, $\sigma ((\T _i)_{l(\T _i)}) \leq \sigma (\T _i)<\alpha  _{\bb{A}_2}(\T _{i+1})\leq \alpha  _{\bb{A}_3}((\T _{i+1})_{l(\T _{i+1})})$, where the last inequality is by Lemma ~\ref{alphaless}.  Thus, $((\T _i )_{l(\T _i)},(\T _{i+1})_{l(\T _{i+1})})$ is an $\LL{\bb{A}_2}$-vector.  Furthermore,
\begin{eqnarray*}
\sigma ((\T _i)_{l(\T _i)-j})&\leq &\sigma (\T _i)-j\mbox{ by Lemma \ref{sigmaless}}\\
                                             &<&\alpha _{\A _2}(\T _{i+1})-j\\
                                             &\leq&\alpha _{\A _3}((\T _{i+1})_1)+l(\T _{i+1})-j-1\\
                                             &\leq& \alpha _{\A _3}((\T _{i+1})_{l(\T _{i+1})-j})
\end{eqnarray*}

\noindent
Thus, each $((\T _i )_{l(\T _i)-j},(\T _{i+1})_{l(\T _{i+1})-j})$ is an $lpp_{\leq}(\mathbb{A}_2)$-type vector for $0\leq j\leq l(\T _i)-1$.

Thus, by the induction hypothesis (and since $l(\T _i)\leq l(\T _{i+1})$), 
\begin{eqnarray*}
W_{(\T _i )_{l(\T _i)}}&\supsetneq& W_{(\T _{i+1})_{l(\T _{i+1})})};\\
W_{(\T _i )_{l(\T _i)-1}}&\supsetneq& W_{(\T _{i+1})_{l(\T _{i+1})-1})};\\ 
&\vdots & \\
W_{(\T _i )_1}&\supsetneq& W_{(\T _{i+1})_{l(\T _{i+1}) -l(\T _i)+1})}.
\end{eqnarray*}
Thus, 
\begin{eqnarray*}
W_{\T _{i+1} } &:=&\langle x_2 ^{l(\T _{i+1})}, x_2 ^{l(\T _{i+1})-1}(W_{(\T _{i+1})_1}), \ldots ,x_2^{l(\T _i)}W_{(\T_{i+1})_{l(\T _{i+1})-l(\T _i)}},\\
               & &x_2^{l(\T _i) -1}W_{(\T _{i+1})_{l(\T _{i+1}) -l(\T _i)+1}},\ldots ,W_{(\T _{i+1})_{l(\T _{i+1})})}\rangle \\
               &\subsetneq& \langle x_2^{l(\T _i)},x_2^{l(\T _i )-1} W_{(\T _i)_1}, \ldots ,W_{(\T _i)_{l(\T _i)}}\rangle\\
               &=&W_{\T _i}
\end{eqnarray*}
\end{proof}

\begin{theorem} \label{DeltaH=codim1} Let $\T$ be an $\LL{\bb{A}}$-vector.  Then $H(R/W_{\T}) = H_{\T }$.
\end{theorem}

\noindent
\begin{proof} We use induction on $n$, the length of $\bb{A}$.  If $n=1$, the result is clear.  So suppose that $n>1$.  Let $\T =(\T _1, \ldots ,\T _s)$.  Let $R=k[x_1,\ldots ,x_n]$.  Then $W_{\T} =\langle x_1^s, x_1^{s-1}\overline{W_{\T _1}}, \ldots ,\overline{W_{\T _s}}\rangle $.  It is enough to show that 
  $$\mbox{codim } (W_{\T} )_d= \sum_{e=1}^s \mbox{ codim
  }(\overline{W_{\T _e}})_{d-s+e}.$$
  Now,
$$\mbox{codim }(W_{\T} )_d =\# \{ \mbox{monomials in $R_d$ not in }W _{\T} \}.$$  Let $M$ be the set of all monomials of $R$ not in $W_{\T}$, and let $T=k[x_2,\ldots ,x_n]$.  Then,
\begin{eqnarray*}
M&\subseteq&
\hspace{.35cm}\{ \mbox{monomials in } T \mbox{ not in }\overline{W _{\T_s}}\}\\
 & &\stackrel{\cdot}{\cup} \{ x_1\cdot (\mbox{monomials in }T\mbox{ not in }\overline{W _{\T_{s-1}}})\}\\
 & &\stackrel{\cdot}{\cup} \hspace{1in}\cdots \\
 & &\stackrel{\cdot}{\cup} \{ x_1^{s-1} \cdot (\mbox{monomials in }T\mbox{ not in }\overline{W _{\T_1}})\}
\end{eqnarray*}
We will show equality.  Certainly, any monomial of $T$ that is not in
$\overline{W_{\T _s}}$ cannot be in $W_{\T} $.  Consider any monomial
$m$ of $x_1 ^{s-i} T$ that is not in $x_1 ^{s-i}\overline{W_{\T _i}}$.
By Lemma ~\ref{containment}, $\overline{W_{\T _j}}\subseteq
\overline{W_{\T _i}}$ for all $j\geq i$.  Write $m=x_1^{s-i}p$, where
$p\in k[x_2,\ldots ,x_n]$.  Now, if we had $m\in W_{\T}$, then we
would have $\frac{m}{x_1^{j-i}}\in x_1^{s-j}\overline{W_{\T _j}}$ for
some $j>i$.  In other words, $m=x_1^{s-i}p$ for some $p\in
\overline{W_{\T _j}}$, and some $j>i$.  This contradicts that
$\overline{W_{\T _j}}\subseteq \overline{W_{\T _i}}$ for all $j\geq
i$.
\end{proof}

So far, we have seen that if $\T$ is an $\LL{\bb{A}}$-vector, then
$W_{\T}$ is an $\LL{\bb{A}}$-ideal with $H(R/W_{\T})=H_{\T}$, $\alpha
_{\mathbb{A}}(H)=\alpha_{\mathbb{A}}(\T )$ and $\sigma (H)=\sigma (\T
)$.  In particular, $H(R/W_{\T})$ is an $\LL{\bb{A}}$-sequence.  We
now wish to show that given any $\LL{\bb{A}}$-sequence $H$, we can
obtain an $\LL{\bb{A}}$-vector $\T$, and furthermore that the function
$H\rightarrow \T$ and the function $\T \rightarrow H_{\T}$ are
inverses of each other.

We begin by decomposing a given $\LL{\bb{A}}$-sequence $S$ into two
``smaller'' such sequences $S_1$ and $S_1'$ by using a decomposition
similar to that used by Geramita, Maroscia and Roberts in \cite{gmr}.
Suppose $S=1\ b_1\ b_2\ b_3 \ldots $, where $b_1 \geq 2$.

Put $e_i =\binom{a_n-1,a_{n-1}-1,\ldots ,a_{n-(b_1-2)}-1}{i}$ and $c_i
=b_{i+1}-e_{i+1}$.  Define $S_1$ as follows:
\begin{enumerate}
\item if $c_i\geq 0$ for all $i$, set $S_1(i)=c_i$ for all $i$;
\item if $c_i\geq 0$ for all $i\leq h-1$ and $c_h<0$, then set $S_1 =c_0\ c_1\ldots c_{h-1} 0\rightarrow$.
\end{enumerate}

In any case, we let $h$ (possibly infinite) be the smallest integer
for which $c_h <0$.  Then define $S_1'$ as follows:

$$S_1'(i) = \left\{ \begin{array}{lll}
                     e_i&\mbox{if}&i\leq h\\
                     b_i&\mbox{if}&i\geq h+1.
                     \end{array} \right.
$$

From the definition of $S_1$ and $S_1'$, it is clear that
$S(i)=S_1'(i) + S_1(i-1)$.

\begin{theorem} Let $S =\{ b_i\}_{i\geq 0}$ be an $\LL{\bb{A}}$-sequence.  Let $S_1$ and $S_1'$ be constructed as above.  Then $S_1$ and $S_1'$ are $\LL{\bb{A}}$-sequences.
\end{theorem}

\noindent
\begin{proof}  Using the Macaulayesque notation for the generalized binomial coefficients, the proof of this statement follows word for word the proof of \cite[Theorem 3.2]{gmr}, so we omit it.
\end{proof}

Before showing the correspondence between $\LL{\bb{A}}$-vectors and Hilbert functions of $\LL{\bb{A}}$-ideals, we need the following lemma.

\begin{lemma} \label{comparealphas} Let $\A =\{ a_1,\ldots ,a_n\}$ and let $S$ be an $\LL{\A}$-sequence, and $S_1$ obtained from $S$ as above.  Suppose that $S(1)=n$.  Then $\alpha _{\A}(S_1)<\alpha _{\A}(S)$.
\end{lemma}

\begin{proof} If $S_1(1)<S(1)$, then $\alpha _{\mathbb{A}}(S_1)=1<\alpha _{\mathbb{A}}(S)$, so suppose that $S_1(1)=S(1)$.  We consider three cases.

\noindent
{\bf Case 1:} $\alpha_{\mathbb{A}}(S)\leq h$. We again use the notation that $a'_i=a_i-1$. Then 
\begin{eqnarray*} 
S_1(\alpha _{\A}(S)-1)&=&b_{\alpha _{\A}(S)}-e_{\alpha _{\A}(S)}\\ 
             &<&\binom{a'_n,a'_{n-1},\ldots ,a'_1}{\alpha _{\A}(S)} - \binom{a'_n,a'_{n-1},\ldots ,a'_2}{\alpha _{\A}(S)}\\ 
              &=&  \binom{a'_n,a'_{n-1},\ldots ,a'_2}{\alpha _{\A}(S)-1}  +\binom{a'_n,a'_{n-1},\ldots ,a'_2}{\alpha _{\A}(S)-2} +\cdots + \binom{a'_n,a'_{n-1},\ldots ,a'_2}{\alpha _{\A}(S)-a'_1} \\
              &\leq &  \binom{a'_n,a'_{n-1},\ldots ,a'_2}{\alpha _{\A}(S)-1}  +\cdots +\binom{a'_n,a'_{n-1},\ldots ,a'_2}{\alpha _{\A}(S)-a'_1} + \binom{a'_n,a'_{n-1},\ldots ,a'_2}{\alpha _{\A}(S)-1-a'_1} \\
&=&\binom{a'_n,a'_{n-1},\ldots ,a'_1}{\alpha _{\A}(S)-1}
\end{eqnarray*}
So, $\alpha _{\A}(S_1)\leq \alpha _{\A}(S)-1$.

\noindent
{\bf Case 2:} $h+1\leq \alpha _{\A}(S)<\infty$.  Then
$S_1(\alpha_{\A}(S)-1)=0<\binom{a'_n,a'_{n-1},\ldots ,a'_1}{\alpha
  _{\A}(S)-1}$, so $\alpha_{\A}(S_1)<\alpha_{\A}(S)$.

\noindent
{\bf Case 3:} $\alpha _{\mathbb{A}}(S)=\infty$.  Then
$S(i)=b_i=\binom{a'_n,\ldots ,a'_1}{i}$ and in particular, $b_1 =
\binom{a'_n,\ldots ,a'_1}{1}= n$, so $e_i = \binom{a'_n,\ldots
  ,a'_2}{i}$. Then,
\begin{eqnarray*}
S_1(i) &=& b_{i+1}-e_{i+1}\\
            &=&\binom{a'_n,\ldots ,a'_1}{i+1} - \binom{a'_n,\ldots ,a'_2}{i+1}\\
            &=&\binom{a'_n,\ldots ,a'_2}{i}+\cdots +\binom{a'_n,\ldots ,a'_2}{i-a'_1}\\
            &=&\binom{a'_n,\ldots ,a'_2,a'_1 -1}{i} 
\end{eqnarray*}

and hence $\alpha_{\A}(S_1)<\infty =\alpha _{\A}(S)$.
\end{proof}

\begin{theorem} 
  There is a 1-1 correspondence between $\LL{\bb{A}}$-vectors and
  Hilbert functions of $\LL{\bb{A}}$-ideals, where if $\T$ corresponds
  to $H$ (we write $\T\leftrightarrow H$), then $\alpha _{\mathbb{A}}
  (\T) =\alpha _{\mathbb{A}} (H )$ and $\sigma (\T) =\sigma (H)$.
\end{theorem}

\begin{proof}
  We first show that the map $\T \rightarrow H_{\T}$ is 1-1.  We
  already know that it preserves $\sigma$ and $\alpha _{\mathbb{A}}$
  and that it does map $\LL{\bb{A}}$-vectors to
  $\LL{\bb{A}}$-sequences.  We use induction on $n$, the base case
  $n=1$ being trivial.
  
  So suppose that $\T \rightarrow H$ and $\T ' \rightarrow H$. We
  first reduce to the case where $\T _{l(\T )}\neq \T
  _{c.i.(a_2,\ldots ,a_n)}$ and $\T ' _{l(\T ')}\neq \T
  _{c.i.(a_2,\ldots ,a_n)}$.

Suppose that $\T =(\T _1, \ldots ,\T _u,\T _u, \ldots ,\T _u)$ and $\T
'=(\T _1 ', \ldots ,\T _v ',\T _v ', \ldots ,\T _v ')$ where $\T _u
=\T _{c.i.(\A _2)}$.  Then $\sigma (\T )=\sigma (c.i.(a_2,\ldots
,a_n))+\# \T _u$'s -1.

If $\T '_v\neq \T _{c.i.(\A _2)}$, then $\sigma (\T ')=\sigma (\T
'_v)<\sigma (c.i.(a_2,\ldots ,a_n))$, contradicting that $\sigma (\T
')=\sigma (\T )$.  So, $\T ' _v =\T _{c.i.(a\A _2)}$ and $\sigma (\T
')=\sigma (c.i.(a_2,\ldots ,a_n))+\# \T _v'$'s $-1$.  Then $\# \T
_u$'s $=\# \T '_v$'s.  So, we also have $(\T _1, \ldots ,\T _{u-1})$
and $(\T ' _1, \ldots ,\T '_{v-1})$ get mapped to the same Hilbert
function.  Thus, we may assume that $\T _{l(\T )}\neq \T _{c.i.(\A
  _2)}$ and $\T ' _{l(\T ')}\neq \T _{c.i.(\A _2)}$.

So let $\T =(\T _1, \ldots ,\T _u )$ and $\T ' =(\T ' _1, \ldots ,\T
'_v )$. Since $\alpha _{\mathbb{A}} (\T )=\alpha _{\mathbb{A}}
(H)=\alpha _{\mathbb{A}} (\T ')$, we have $u=v$.  From here, the
argument that $\T =\T '$ follows word for word the argument in
\cite[Theorem 2.6]{ghs}, so we omit it.

Now, we define the map $H\rightarrow \T$ inductively as follows:

If $n=1$, then $H=1\ 1\ldots 1\ 0\rightarrow$ where there are $d$ 1's,
for some $d\leq a_1$.  So put $H\rightarrow \T =(d)$.

If $n>1$, we may as well assume that $a_i\geq 2$ for all $i$, and that
$H(1)=n$, for if $H(1)<n$, then we claim that $H$ is also an $\LL{\A
  _2}$-sequence.  To see this, consider the Macaulayesque rectangle
used to construct $\LL{\A}$-sequences, where the $i^{\tiny {\rm th}}$
row consists of $\binom{a_n-1,\ldots ,a_{n-i+1}-1}{j}$ for $j\geq 0$.
So if $H(1)\leq n-1$, then $H$ cannot in any degree occur below the
row consisting of $\binom{a_n-1,\ldots ,a_2-1}{j}$ for $j\geq 0$.  So
$H$ is also an $\LL{\A _2}$-sequence.  Thus, in this case we may use
induction on $n$.

Now, decompose $H$ into $H_1$ and $H'_1$.  By induction on $n$, send
$H'_1\rightarrow \T '_1$.  By Lemma ~\ref{comparealphas}, $\alpha
_{\bb{A}} (H_1)<\alpha _{\bb{A}}(H)$, so by induction on $\alpha
_{\bb{A}}$ (the base case $\alpha _{\bb{A}} =1$ being the induction
hypothesis on $n$), we send $H_1\rightarrow \T _1=((\T _1)_1,\ldots
,(\T _1)_{l(\T _1)})$.  Then send $H\rightarrow ((\T _1)_1,\ldots ,(\T
_1)_{l(\T _1)}, \T '_1)$.  This is an $\LL{\A}$-vector, since
\begin{eqnarray*}
\sigma ((\T_1)_{l(\T _1)})\leq \sigma (\T _1)&=&\sigma (H_1)\hspace{.25cm} \mbox{by induction}\\
                                                                             &\leq&h \hspace{.25cm}\mbox{by construction of $H_1$}\\
                                                                             &<&\alpha _{\A _2}(H'_1)\hspace{.25cm}\mbox{by construction of $H'_1$}.
\end{eqnarray*}

Next we claim that $H\rightarrow \T \rightarrow H$ is the identity
map.  This is clearly true when $n=1$, so we use induction on $n$ and
assume that $n>1$. Note that if $H\rightarrow \T =(\T_1,\ldots
,\T_u)$, we must have $H_1\rightarrow (\T_1,\ldots , \T_{u-1})$ and
$H_1'\rightarrow \T_u$, by definition.  Then

\begin{eqnarray*}
H\rightarrow \T=(\T_1,\ldots ,\T_u) &\rightarrow&H_{\T_u}(i) +H_{(\T_1,\ldots ,\T_{u-1})}(i-1)\hspace{.25cm}\mbox{by definition}.\\
         &=&H'_1(i)+H_1(i-1)\hspace{.25cm}\mbox{by induction since}\\
         & &\hspace{.25cm} H'_1\rightarrow \T _u\mbox{ and }H_1 \rightarrow (\T_1,\ldots , \T_{u-1})\\
         &=&H(i).
\end{eqnarray*}

This, together with $\T \rightarrow H_\T$ being 1-1 shows that $\T
\rightarrow H_{\T}$ and $H\rightarrow \T$ are inverses of each other.
\end{proof}

\begin{example} Let 
  $$\T =(\T _1,\T _2,\T _3,\T _4)
  =((1,2),(1,3,4),(2,3,6,6),(5,6,6,6))$$
  be an $lpp_{\leq}(\{ 4,4,6\}
  )$-vector.  Then letting $\T \rightarrow H$ and $\T _i\rightarrow
  H_i$, we have

$$
\begin{array}{ccccccccccc}
H_4:&1&2&3&4&4&4&3&2&0&\rightarrow\\
H_3:& &1&2&3&4&4&2&1&0&\rightarrow\\
H_2:& & &1&2&3&2&0&\rightarrow& & \\
H_1:& & & &1&2&0&\rightarrow& & & \\
\hline\\
H:  &1&3&6&10&13&10&5&3&0&\rightarrow
\end{array}
$$

Now, beginning with $H=$ 1 3 6 10 13 10 5 3 0 $\rightarrow$, an
$\LL{\bb{A}}$-sequence, we have:

$$
\begin{array}{cccccccccccc}
b_i:&1&3&6&10&13&10&5&3&0&\rightarrow& \\
e_i:&1&2&3& 4& 4& 4&3&2&1&0&\rightarrow\\
\hline\\
c_i:& &1&3& 6& 9& 6&2&1&-1&0&\rightarrow
\end{array}
$$

So, $S_1 =$1 3 6 9 6 2 1 0 $\rightarrow$ and $S'_1 =$1 2 3 4 4 4 3 2 0
$\rightarrow$.

Continuing, we decompose $S$ as 

$$
\begin{array}{cccccccccc}
1&2&3&4&4&4&3&2&0&\rightarrow\\
 &1&2&3&4&4&2&1&0&\rightarrow\\
 & &1&2&3&2&0&\rightarrow& & \\
 & & &1&2&0&\rightarrow& & & \\
\end{array}
$$

We decompose each of these further to obtain:
$$
\begin{array}{lcl}
1\ 2\ 3\ 4\ 4\ 4\ 3\ 2\ 0\ \rightarrow&\longleftrightarrow&(5,6,6,6)\\
1\ 2\ 3\ 4\ 4\ 2\ 1\ 0\ \rightarrow&\longleftrightarrow&(2,3,6,6)\\
1\ 2\ 3\ 2\ 0\ \rightarrow&\longleftrightarrow&(1,3,4)\\
1\ 2\ 0\ \rightarrow&\longleftrightarrow&(1,2).
\end{array}
$$

So we indeed obtain $\T$ back from $H$.
\end{example}

\section{ideal colon}
In this section, our goal is to show that the residual of an
$\LL{\A}$-ideal in the complete intersection of type $(a_1, \ldots
,a_n)$ is again an $\LL{\bb{A}}$-ideal.

In two variables, where $\A =\{ a,b\}$, the residual of an $LPP(\A
)$-ideal inside the $c.i.(a,b)$ is always a lex ideal, namely
$$
\langle x^a, y^b\rangle : \langle x^a,
x^{a-1}y^{d_1},x^{a-2}y^{d_2}, \ldots ,x^{a-s}y^{d_s}, y^b\rangle $$
$$=\langle x^s,x^{s-1}y^{b-d_s},\ldots ,xy^{b-d_2},y^{b-d_1}\rangle
.$$

As before, we associate to the $LPP(\A )$-ideal 
$$\langle x^a, x^{a-1}y^{d_1},x^{a-2}y^{d_2}, \ldots ,x^{a-s}y^{d_s},
y^b\rangle$$
the $\LL{\A }$-vector $\T = (d_1, \ldots ,d_s,b, , \ldots
,b)$, where there are $a-s$ $b$'s, so that the length of $\T$ is $a$.
Then we associate to the residual lex ideal the 2-type vector $(b-d_s,
\ldots ,b-d_1)$.  We can use monomial lifting (see \cite[Theorem
2.2]{ggr}) to associate a finite set of points to each of these
ideals.  The set of points obtained from the lex ideal in this way is
an example of a $k$-configuration.  From the lpp ideal, we obtain the
complement of the $k$-configuration in the $c.i.(a,b)$; this
complementary set of points is an example of a weak $k$-configuration,
as defined in \cite[Definition 2.8]{gps}.  In fact, $\LL{\{ a,b\}
}$-vectors are exactly the ``types'' of weak $k$-configurations that
occur in theorem 2.10 of their paper.  It was this fact that motivated
the definition of $\LL{\{ a,b\} }$-vectors and the generalization to
larger numbers of variables.

\begin{example} 
  The following ideal is LPP(5,7): $I=\langle x^5,
  x^4y,x^3y^3,x^2y^4,y^7\rangle$.  We associate to $I$ the $\LL{\{
    5,7\} }$-vector (1,3,4,7,7).  Then, inside a $c.i.(5,7)$, we draw
  a weak $k$-configuration of type (1,3,4,7,7):
  
  $$
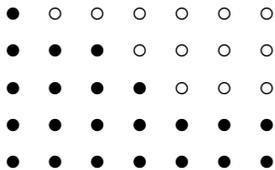
\begin{array}{ccccccc}
\bullet&\circ&\circ&\circ&\circ&\circ&\circ\\
\bullet&\bullet&\bullet&\circ&\circ&\circ&\circ\\
\bullet&\bullet&\bullet&\bullet&\circ&\circ&\circ\\
\bullet&\bullet&\bullet&\bullet&\bullet&\bullet&\bullet\\
\bullet&\bullet&\bullet&\bullet&\bullet&\bullet&\bullet
\end{array}
$$

\noindent
In this case, the complement of the weak $k$-configuration is a
$k$-configuration of type (3,4,6).
\end{example}

The fact that the residual of an $LPP\{ a,b\} $-ideal in the
$c.i.(a,b)$ is a lex ideal provides a proof of the LPP conjecture in
two variables (for another proof, see \cite[Theorems 5.1 and
5.2]{Richert}).  Since $x^a,y^b$ are never minimal generators of the
residual lex ideal, the resolution of the lex plus powers ideal
obtained from dualizing the minimal free resolution of the lex ideal
is in fact minimal (see page 154 of \cite{Migliore}).  Hence, since
lex ideals have extremal resolutions, it follows that the lex plus
powers ideals have extremal resolutions among all ideals containing an
$\{ a,b\}$-regular sequence.

In more than two variables, this argument does not work for several
reasons.  Firstly, the generators of the complete intersection might
be generators of the residual ideal; secondly, even if they were not,
we would not be guaranteed that the resolution obtained by dualizing
was minimal and thirdly, the residual of an $LPP(\A )$-ideal is no
longer necessarily a lex ideal.

In this section, however, we show that the residual of an $\LL{\A
}$-ideal is necessarily another $\LL{\A }$-ideal.  Given an $\LL{\A
}$-vector $\T$, we want to define a residual $\LL{\A }$-vector $\T
^*$.

We also associate to $W_{\T}$, and hence to $\T$, a natural set of
points $\pts$ in $\proj ^n$ contained in a complete intersection of
type $(a_1, \ldots ,a_n)$ obtained from lifting the monomial ideal
$W_{\T}$. Note that we do not need to know that $W_{\T}$ is lex plus
powers in order to associate the set $\pts$ of points in $\proj ^n$;
we only need that it is monomial. Then consider $\pts ^c$, the
complement of $\pts$ in $c.i.(a_1, \ldots ,a_n)$.  We want to define a
dual vector $\T^*$ so that $W_{\T^*}\leftrightarrow \pts ^c$.

\begin{definition}
  If $n=1$, so that $\T$ is an $\LL{\bb{A}}$-vector $(d)$, $d\leq
  a_1$, then $\T ^* := (a_1-d)$ if $d<a_1$; otherwise, we define
  $\T^*=\emptyset$.

If $\T$ is an $\LL{\bb{A}}$-vector $(\T _1, \ldots ,\T _u)$ and if $u<
a_1$, then
$$\T ^* := ((\T _u)^*, \ldots ,(\T _1)^*, \T _{c.i.(\A _2)}, \ldots ,
\T _{c.i.(\A _2)}),$$
where there are $a_1 -u$ $\T_{c.i.(\A _2)}$'s;
otherwise, $\T ^* =((\T_u)^*, \ldots ,(\T _1)^*)$. In particular,
$(\T_u)^* =(\T^*)_1$ unless $(\T _u)^* = \emptyset$.

We also define $W_{\T ^*}$ in the same way we defined $W_{\T}$. While
we do not yet know that $W_{\T^*}$, defined in this way, is lex plus
powers, we do know that it is a monomial ideal and so we can associate
a set of points to $\T^*$ by monomial lifting.
\end{definition}

\begin{remark} 
  With this definition of $\T ^*$, we see that if $\pts
  \leftrightarrow W_{\T}$, we indeed have $\pts ^c \leftrightarrow
  W_{\T^*}$.  Note that we can define $l(\T^*)$, $\alpha _{\A}(\T ^*)$
  and $\sigma (\T ^*)$, just as we defined these parameters for $\T$,
  even before knowing that $\T ^*$ is an $\LL{\A}$-vector; we also put
  $\alpha (\emptyset )=\sigma (\emptyset )=0$.  Furthermore, if we
  perform the same operation on $\T ^*$ as we did on $\T$ to obtain
  $\T ^*$, we obtain $\T$ back.  We write this as $(\T ^*)^*=\T$.  As
  well, it is clear that $l(\T )<a_1 \Leftrightarrow (\T ^*)_{l(\T
    ^*)} = \T_{c.i.(\A _2)}$.
\end{remark}

We want to show that if $\T$ is an $\LL{\bb{A}}$-vector, then so is
$\T^*$.

{\bf Notation:} In what follows, we remove the subscript $\bb{A}$ from
the $\alpha$ notation and assume it to be understood.  So we write
$\alpha (\T )$ for $\alpha _{\bb{A}}(\T )$, $\alpha (\T _i)$ for
$\alpha _{\bb{A}_2}(\T _i)$, $\alpha ((\T _i)_j)$ for $\alpha
_{\bb{A}_3}((\T _i)_j)$, etc., assuming the subscript is understood.

\begin{lemma} \label{a+s} 
  Let $\T$ be an $\LL{\A }$-vector.  Then
  $$\alpha (\T )+\sigma (\T ^*)=\sigma (c.i.(a_1, \ldots ,a_n))
  =\sigma (\T ) +\alpha (\T ^*).$$
\end{lemma}

\begin{proof}
  When $n=1$, the result is trivial.  To show that $\alpha (\T
  )+\sigma (\T^*)=\sigma (c.i.(\mathbb{A}))$, we consider the two
  cases $l(\T )<a_1$ and $l(\T )=a_1$.
  
  If $l(\T ) <a_1$, then $(\T ^*)_{l(\T ^*)} = \T _{c.i.(a_2, \ldots
    ,a_n)}$.  Then, $\alpha (\T ) +\sigma (\T ^*) = l(\T ) +\sigma
  (c.i.(a_2, \ldots ,a_n)) +s-1$ where $s = \# \{ i | (\T ^*)_i =\T
  _{c.i.(\A _2)}\}$.  But $l(\T ) +s =a_1$, so $\alpha (\T ) + \sigma
  (\T ^*) = a_1 +\sigma (c.i.(a_2, \ldots ,a_n))-1 = \sigma (c.i.(a_1,
  \ldots ,a_n))$.

If $l(\T ) =a_1$, then $(\T ^*)_{l(\T ^*)}\neq \T _{c.i.(\A _2)}$.  So
$\alpha (\T )+\sigma (\T ^*) = a_1+\alpha (\T _1) -1 +\sigma ((\T
^*)_{l(\T ^*)}).$ But $(\T ^*)_{l(\T ^*)} =(\T _1)^*$ so $\alpha (\T
_1) + \sigma ((\T ^*)_{l(\T ^*)}) =\sigma (c.i.(a_2, \ldots ,a_n))$ by
induction, so $\alpha(\T ) +\sigma (\T ^*)=a_1 -1 +\sigma (c.i.(a_2,
\ldots ,a_n))=\sigma (c.i.(a_1, \ldots ,a_n))$

To show that $\sigma (\T ) +\alpha (\T ^*) = \sigma (c.i.(a_1, \ldots
,a_n))$, we consider the two cases $\T_{a_1} = \T_{c.i.(a_2,\ldots
  ,a_n)}$ and $\T_{a_1} \neq \T_{c.i.(a_2,\ldots ,a_n)}$.

If $\T_{a_1} = \T_{c.i.(a_2,\ldots ,a_n)}$, then $l(\T^*) <a_1$, so
$\alpha (\T^*) =l(\T^*)$. Furthermore, $\sigma (\T) =\sigma
(c.i.(\mathbb{A}_2)) +s-1$, where $s$ is the number of integers $i$
such that $\T_i =\T_{c.i.(\mathbb{A}_2)}$. But $l(\T^*)+s = a_1$, so
$\alpha (\T^*)+\sigma (\T) = a_1 +\sigma (c.i.(\mathbb{A}_2)) -1 =
\sigma (c.i.(\mathbb{A})).$

If $\T_{a_1} \neq \T_{c.i.(a_2,\ldots ,a_n)}$, then $\l(\T^*) = a_1$.
So $\alpha (\T^*) = a_1+\alpha _{\mathbb{A}_2}((\T^*)_1)-1 = a_1
+\alpha _{\mathbb{A}_2}((\T_{l(\T )})^*) -1$, since $(\T_{l(\T
  )})^*\neq \emptyset$. Furthermore, $\sigma (\T) = \sigma
(\T_{l(\T)})$. By the induction hypothesis, $\alpha
_{\mathbb{A}_2}((\T_{l(\T )})^*) +\sigma (\T_{l(\T)}) = \sigma
(c.i.(\mathbb{A}_2))$, so $\sigma (\T) + \alpha (\T^*) = a_1 +\sigma
(c.i.(\mathbb{A}_2)-1 = \sigma (c.i.(\mathbb{A})).$

\end{proof}

\begin{theorem} \label{s<adual}
  Let $\mS$ and $\T$ be $\LL{\bb{A}}$-vectors.  Then $\sigma (\mS
  )<\alpha (\T )\Rightarrow \sigma (\T ^*)<\alpha (\mS^*).$
\end{theorem}

\begin{proof} We consider several cases:
\vspace{.25cm}

\noindent
{\bf Case 1:} $\mS_{l(\mS)} =\T_{c.i.(a_2, \ldots ,a_n)}$ and $l(\T
)<a_1$. Then $l(\mS^*)<a_1$ and $(\T^*)_{l(\T^*)} =\T_{c.i.(a_2,\ldots
  ,a_n)}$.  Then since $\sigma (\mS) <\alpha (\T)$, we have $\sigma
(\mS_{l(\mS)}) +s-1<l(\T)$ where $s= \# \{ i |\mS_i =\mS_{l(\mS)}\}$.
Now, $\sigma (\T^*)<\alpha (\mS^*)\Leftrightarrow \sigma
((\T^*)_{l(\T^*)}) +t-1<l(\mS^*)$ where $t=\# \{ i |(\T^*)_i =
(\T^*)_{l(\T^*)}\}$. But $\sigma (\mS_{l(\mS)}) =\sigma
(\T^*_{l(\T^*)})=\sigma (c.i.(a_2, \ldots ,a_n))$, so it is enough to
show that $l(\T)-s =l(\mS^*)-t$.  But, $l(\T )+t =a_1 =l(\mS^*)+s$, so
we are done in this case.  \vspace{.25cm}

\noindent
{\bf Case 2:} $\mS_{l(\mS)} =\T _{c.i.(a_2, \ldots ,a_n)}$ and
$l(\T)=a_1$. Then $l(\mS^*)<a_1$ and $(\T^*)_{l(\T^*)}\neq \T
_{c.i.(a_2, \ldots ,a_n)}$.  So, $\sigma (\T^*)=\sigma
((\T^*)_{l(\T^*)})$ and $\alpha (\mS^*)=l(\mS^*)$.  So we need to show
that $\sigma ((\T^*)_{l(\T^*)})<l(\mS^*)$.  Note that $l(\mS^*)=a_1
-s$ where $s =\# \{ i |\mS_i =\mS_{l(\mS)}\} $. Since $\sigma
(\mS)<\alpha (\T)$, we have $\sigma (\mS_{l(\mS)}) +s-1 <l(\T) +\alpha
(\T_1) -1$.  Rewriting this, we obtain $\sigma (\mS_{l(\mS)}) - \alpha
(\T_1)<l(\T) -s = a_1 -s.$ But $(\T_1)^* = (\T^*)_{l(\T^*)}$ since
$l(\T)=a_1$, so by Lemma ~\ref{a+s}, $\sigma ((\T^*)_{l(\T^*)}) =
\sigma (c.i.(\mathbb{A}_2)) - \alpha (\T_1) <a_1-s=l(S^*)$, as
required.  \vspace{.25cm}

\noindent
{\bf Case 3:} $\mS_{l(\mS)}\neq \T_{c.i.(a_2, \ldots ,a_n)}$ and
$l(\T) <a_1$.  Then $l(\mS^*)=a_1$ and $(\T^*)_{l(\T^*)}=\T
_{c.i.(a_2,\ldots ,a_n)}$.  Let $s=\# \{ i |(\T^*)_i
=(\T^*)_{l(\T^*)}\}$.  We need to show that $\sigma
((\T^*)_{l(\T^*)})+s-1<a_1+\alpha ((\mS^*)_1)-1$ or in other words,
$\sigma (c.i.(a_2,\ldots ,a_n)) - \alpha ((\mS^*)_1)<a_1-s =l(\T).$
But $l(\mS^*) =a_1$, so $((\mS^*)_1)^*=\mS_{l(\mS)}$ and hence
$(\mS^*)_1=(\mS_{l(\mS)})^*$.  So by Lemma ~\ref{a+s} applied to
$\mS_{l(\mS)}$, the left hand side of this last inequality is $\sigma
(\mS_{l(\mS )})$. But $\sigma (\mS_{l(\mS)})=\sigma (\mS) <\alpha (\T
) = l(\T )$, so we are done in this case.

\vspace{.25cm}

\noindent
{\bf Case 4:} $\mS_{l(\mS)}\neq \T_{c.i.(a_2, \ldots ,a_n)}$ and
$l(\T)=a_1$. Then, $l(\mS^*)=a_1$ and $(\T^*)_{l(\T^*)}\neq \T
_{c.i.(a_2,\ldots ,a_n)}$. Since $\sigma (\mS)<\alpha (\T)$, we have
$\sigma (\mS_{l(\mS)})<a_1 +\alpha ((\T)_1)-1$. We need to show that
$\sigma (\T^*)<\alpha (\mS^*)$, in other words, $\sigma
((\T^*)_{l(\T^*)})<a_1+\alpha ((\mS^*)_1)-1.$ So it is enough to show
that $\sigma (\mS_{l(\mS)})-\alpha (\T_1)=\sigma
((\T^*)_{l(\T^*)})-\alpha ((\mS^*)_1),$ in other words $\sigma
(\mS_{l(\mS)})+\alpha ((\mS^*)_1)=\sigma ((\T^*)_{l(\T^*)})+\alpha
(\T_1).$ But $l(\mS^*)=a_1$, so $((\mS^*)_1)^*=\mS_{l(\mS)}$ and
$l(\T)=a_1$, so $(\T_1)^*=(\T^*)_{l(\T^*)}$.  Thus by Lemma ~\ref{a+s}
applied to $\mS_{l(\mS)}$ and $\T_1$, $\sigma (\mS_{l(\mS)}) +\alpha
((\mS^*)_1)=\sigma (c.i.(a_2, \ldots ,a_n))=\sigma
((\T^*)_{l(\T^*)})+\alpha (\T_1)$, as required.
\end{proof}

\begin{corollary} \label{c:residuals} 
  If $\T$ is an $\LL{\bb{A}}$-vector, then so is $\T ^*$.  In
  particular, if $I$ is an $\LL{\bb{A}}$-ideal, then so is $\langle
  x_1^{a_1}, \ldots ,x_n^{a_n}\rangle :I$.
\end{corollary}

\begin{remark}
  Chris Francisco has also discovered a (quite different) proof of
  this result.
\end{remark}

\begin{proof}
  Let $\mathbb{A}=\{ a_1,\ldots ,a_n\}.$ If $n=1$, the result is
  obvious, so assume $n>1$ and let $\T = (\T _1, \ldots ,\T_u)$ be an
  $\LL{\mathbb{A}}$-vector so that $u\leq a_1$, $u\leq l(\T _u)$, each
  $\T_i$ is an $\LL{\mathbb{A}_2}$-vector and $\sigma (\T _i)<\alpha
  _{\mathbb{A}_2}(\T _{i+1})$ for $1\leq i\leq u-1$.  Then
  $$\T^* = ((\T _u)^*, \ldots ,(\T_1)^*,
  \T_{c.i.(\mathbb{A}_2)},\ldots , \T_{c.i.(\mathbb{A}_2)})$$
  where
  there are $a_1-u$ (possibly 0) $\T_{c.i.(\mathbb{A}_2)}$'s.  By the
  induction hypothesis, each $\T_i^*$ is an $\LL{\mathbb{A}_2}$-vector
  and $l(\T^*)\leq a_1$ by construction.  To see that $l(\T^*)\leq
  l((\T^*)_{l(\T^*)})$, we consider two cases.
  
  {\bf Case 1:} $u<a_1$.  Then
  $(\T^*)_{l(\T^*)}=\T_{c.i.(\mathbb{A}_2)}$ and $l((\T^*)_{l(\T^*)})
  =a_2 \geq a_1\geq l(\T^*)$.
  
  {\bf Case 2:} $u=a_1$.  Let $\T =(\T _1, \ldots ,\T_u) =(\T_1,
  \ldots ,\T _s,\T_{c.i.(\mathbb{A}_2)},\ldots ,
  \T_{c.i.(\mathbb{A}_2)})$ where $\T _s \neq \T_{c.i.(\mathbb{A}_2)}$
  and $s\leq u = a_1$.  Note that if $\T = \T_{c.i.(\mathbb{A})}$,
  then $s=0$. Then $l(\T^*)=s$ and $(\T^*)_{l(\T ^*)} =(\T_1)^*$, so
  we need to show that $l((\T_1)^*)\geq s.$ First note that since
  $\sigma (\T _1)<\alpha _{\mathbb{A}_2}(\T _2)\leq\sigma (\T
  _2)<\cdots < \alpha _{\mathbb{A}_2}(\T _s)\leq\sigma (\T _s),$ we
  have $\sigma (\T _1)\leq \sigma (\T _s) -s+1\leq
  \sigma(c.i.(\mathbb{A}_2))-1-s+1=\sigma(c.i.(\mathbb{A}_2))-s$.
  
  Let $t$ be the number of $\T_{c.i.(\mathbb{A}_3)}$'s in $\T_1$. Then
  $\sigma (\T _1) = \sigma (c.i.(\mathbb{A}_3)) +t-1 \leq \sigma (c.i.
  (\mathbb{A}_2)) - s = \sigma (c.i.(\mathbb{A}_3)) +a_2-1-s$. So,
  $t\leq a_2-s$; that is, $s\leq a_2 - t =l((\T_1)^*)$, as required.
  
  Thus, it only remains to prove that $\sigma ((\T_i)^*)<\alpha
  ((\T_{i-1})^*)$, but this is the content of Theorem ~\ref{s<adual}.
\end{proof}

\section{Applications of the theorem for colon ideals}\label{apply} 

The fact that the residual of a lex plus powers ideal is again lex
plus powers allows us to prove the (moral) converse to the following
theorem in \cite{Richert}:
\begin{theorem}
  Let $L$ be $\LA$ for some $\A=\a$, and $I$ be an ideal containing an
  $\A$-regular sequence such that $H(R/L)=H(R/I)$. If LPPH holds, then
  $\dim_k(\soc(L)_d)\ge \dim_k(\soc(I)_d)$ for all $d$, where
  $\soc(L)_d$ refers to the $d$th graded piece of the socle of $R/L$
  (and similarly for $I$).
\end{theorem}

We will here demonstrate that if lex plus powers ideals can be shown
to have always largest socles, the LPPH must be true. More precisely,
we will prove that LPPH is equivalent to the following conjecture:

\begin{conjecture}\label{c:socles}
  Suppose that $L$ is $\LA$ for some $\A=\a$, and $I$ is an ideal
  containing an $\A$-regular sequence such that $H(R/L)=H(R/I)$. Then
  $\beta^L_{n,j}\ge \beta^I_{n,j}$ for all $j$.
\end{conjecture}

The proof of the equivalence will require a few lemmas and a
proposition.  We give the following comments to motivate these
preliminary results. Suppose that $L$ is $\LA$ with $\X=\{x_1^{a_1},
\dots , x_n^{a_n}\}$, $I$ contains an $\A$-regular sequence
$\Y=\{y_1^{a_1}, \dots , y_n^{a_n}\}$, and $H(R/L)=H(R/I)$. Our goal
is to compare the socles of $(\X:L)$ and $(\Y:I)$ (via conjecture
\ref{c:socles}) and transfer this comparison to a comparison of the
first graded Betti numbers of $L$ and $I$. By corollary
\ref{c:residuals}, we know that $(\X:L)$ is again a lex plus powers
ideal, so conjecture \ref{c:socles} will apply if we can demonstrate
that $(\Y:I)$ contains a regular sequence in the same degrees as those
of the minimal monomial regular sequence in $(\X:L)$ (note that the
Hilbert functions of the two colon ideals are obviously equal). This
follows from the lemmas below: we first prove (lemma \ref{l:lexseg})
that if $L$ is $\LA$, then the degrees of the minimal monomial regular
sequence in the residual can only drop in degrees for which the colon
consists of a lex segment.  We then use this fact to show (lemma
\ref{l:regsequence}) that $(\Y:I)$ contains a regular sequence in the
degrees of the minimal monomial regular sequence in $(\X:L)$.
Proposition \ref{p:lasttwo} then allows us to compute the first graded
Betti numbers of $L$ and $I$ from the socle degrees of $(\X:L)$ and
$(\Y:I)$ respectively. After these preparations, we will be able to
prove the theorem.

\begin{lemma}\label{l:lexseg}
  Let $L$ be an $\a$-lex plus powers ideal and $\X=\ax$. If $(\X:L)$
  is an $\ap$-lex plus powers ideal with $a'_s<a_s$ for some $1\le s
  \le n$, then $(\X:L)_{a'_s}$ is a lex segment.
\end{lemma}

\begin{proof}
  Note that if $a'_i=a'_s$ for some $i>s$, then $a'_i=a'_s<a_s\le
  a_i$, so we can assume without harm that $a'_s<a'_{s+1}$ or $s=n$.
  It follows that if $m\in(\X:L)_{a'_s}$ and $m<x_s^{a'_s}$, then $m$
  is not a pure power. Because $x_s^{a'_s}$ is a minimal generator,
  $m$ must be a minimal generator as well, and thus it is part of the
  lex segment of $(\X:L)_{a'_s}$.
  
  So it is enough to show that if $m\in R_{a'_s}$ and $m>x_s^{a'_s}$,
  then $m\in (\X:L)$. Note that $s>1$ (otherwise we are finished). If
  $m\not\in(\X:L)$, then there is a minimal monomial generator
  $\lambda\in L$ such that $m\lambda\not\in\X$. It follows that
  $m(i)+\lambda(i)<a_i$ for all $i=1, \dots , n$.  Now, since
  $x_s^{a'_s}\in (\X :L)$, we have $\lambda x_s^{a'_s}\in \langle
  \X\rangle$, and so $\lambda (s) +a'_s\geq a_s$.  In particular, this
  implies that $\lambda(s)>0$.  If $\deg \lambda =d$, then since
  $\lambda (i) <a_i$ for all $i$, $\lambda$ is part of the lex segment
  of $L_d$, and thus if $\lambda'\in R_d$ and $\lambda'>\lambda$, then
  $\lambda'\in L_d$ as well.
  
  Now let $t<s$ be such that $m(t)>0$ (such an element exists because
  $m>x_s^{a'_s}$) and consider the element
  $$\lambda'=x_1^{\lambda(1)+\gamma(1)} \cdots
  x_{s-1}^{\lambda(s-1)+\gamma(s-1)}
  x_{s+1}^{\lambda(s+1)+\gamma(s+1)} \cdots
  x_n^{\lambda(n)+\gamma(n)},$$
  where the $\gamma(i)$ for $i\ne s$ are
  any choice of elements of $\N$ such that $\sum_{i\ne s}
  \gamma(i)=\lambda(s)$, $\gamma(t)\ge 1$, and $\gamma(i)\le m(i)$ for
  all $i\ne s$. Such a choice of $\gamma(i)$ is possible unless
  $\lambda(s)= \sum_{i\ne s} \gamma(i) > \sum_{i\ne s} m(i)=
  \deg(m)-m(s)=a'_s-m(s)$ in which case $\lambda(s)+m(s)>a'_s$, a
  contradiction. The existence of such a $\gamma$, however, also gives
  a contradiction. Because $\lambda'>\lambda$, we have that
  $x_s^{a'_s}\lambda'\in\X$.  But $\lambda '(s)=0$ and $a'_s<a_s$, so
  for some $i\ne s$, $a_i\le \lambda(i)+\gamma(i)\le \lambda(i)+m(i)$.
\end{proof}

\begin{lemma}\label{l:regsequence}
  Suppose that $I$ minimally contains an $\A=\a$-regular sequence
  $\Y$, $H(R/I)$ is $\A$-lpp valid, and let $L$ be the $\A$-lex plus
  powers ideal such that $H(R/I)=H(R/L)$. If $(\X:L)$ is $\ap$-lex
  plus powers, then $(\Y:I)$ contains an $\ap$-regular sequence.
\end{lemma}
\begin{proof}
  Let $t$ be the smallest integer such that $(\Y:I)$ fails to contain
  an $\{a'_1, \dots, a'_t\}$-regular sequence. Thus there is a
  $\b$-regular sequence in $(\X:I)$ such that $b_i\le a'_i$ for $1\le
  i <t$, and $a'_t<b_t\le a_t$. We can choose $\b$ such that $b_t$
  satisfies the second inequality because $(\Y:I)$ contains an
  $\a$-regular sequence by construction and thus certainly contains an
  $\{a'_1, \dots , a'_{t-1}, a_t, a_{t+1}, \dots, a_n\}$-regular
  sequence.  By lemma \ref{l:lexseg}, $a'_t<a_t$ implies that
  $(\X:L)_{a'_t}$ is a lex segment. Consider then the ideals
  $(\X:L)_{a_t'}+\m{x_1, \dots , x_n}^{a_t'+1}$ and
  $(\Y:I)_{a_t'}+\m{x_1,\dots , x_n}^{a_t'+1}$. Both of these ideals
  attain the same Hilbert function, and the former is a lex ideal
  containing a regular sequence of length at least $t$ in degree
  $a_t'$.  It is not difficult to show (see for example, Corollary
  2.13 in \cite{Richert}) that all ideals attaining a given Hilbert
  function contain a regular sequence in the degrees of the minimal
  monomial regular sequence in the lex ideal with that Hilbert
  function. Thus, $(\Y:I)_{a_t'}+\m{x_1,\dots ,x_n}^{a_t'+1}$ must
  also contain a regular sequence of length at least $t$ by degree
  $a_t'$, that is, $(\Y:I)$ must contain a regular sequence in degrees
  $a_1, \dots , a_t$, a contradiction.
\end{proof}

\begin{proposition}\label{p:lasttwo}
  Let $\A=\a$ be a list of degrees, write $|j|$ to denote the number
  of elements of $\A$ equal to $j$, and suppose that $\Y$ is an
  $\A$-regular sequence in an ideal $I\subset R$. Then for all $j$
  there exist $0\le t_j\le |j|$ such that
  $\beta_{n,\omega-j}^{(\Y:I)}=\beta_{1,j}^I-t_j$. Furthermore, if
  $\Y$ is minimally contained in $I$, then $t_j=|j|$ for all $j$.
\end{proposition}
\begin{proof}
  We suppose first that $\Y$ is minimally contained in $I$.  Let
  $\F_{\bullet}$ be a minimal free resolution of $R/I$
  $$\F_{\bullet} := 0 \rightarrow \sum_{j} R^{\beta^I_{n,j}}[-j]
  \xrightarrow[]{\delta_n} \cdots \xrightarrow[]{\delta_{2}} \sum_{j}
  R^{\beta^I_{1,j}}[-j] \xrightarrow[]{\delta_{1}} R \rightarrow 0,$$
  and $\K_{\bullet}$ be the Koszul complex
  $$\K_{\bullet}= 0 \to R[-\omega] \xrightarrow[]{\partial_n} \sum_j
  R^{\beta^K_{n-1,j}}[-j] \xrightarrow[]{\partial_{n-1}} \cdots
  \xrightarrow[]{\partial_{2}} \sum_j R^{\beta^K_{1,j}}[-j]
  \xrightarrow[]{\partial_1} R \to 0$$
  resolving $R/\Y$, where the
  $\beta^K_{i,j}$ are the Betti numbers of the Koszul complex
  resolving $R/\Y$ and $\omega=\sum a_i$. Note that
  $|j|=\beta_{1,j}^K$. The map $\phi:R/\Y\to R/I$ induces a chain map
  {\footnotesize
  $$
\begin{CD} 
  0 @>>> \sum_{j} R^{\beta^I_{n,j}}[-j] @>{\delta_{n}}>> \sum_{j}
  R^{\beta^I_{n-1,j}}[-j] @ >{\delta_{n-1}}>> \cdots @ >{\delta_{2}}>>
  \sum_{j} R^{\beta^I_{1,j}}[-j] @>
  {\delta_{1}}>> R @>>> 0\\
  & & @AA{\phi_n}A  @AA{\phi_{n-1}}A && @AA{\phi_{1}}A @AA{\phi_{0}}A\\
  0 @>>> R[-\omega] @>{\partial_n}>> \sum_j R^{\beta^K_{n-1,j}}[-j] @
  >{\partial_{n-1}}>> \cdots @ >{\partial_2}>>
  \sum_j R^{\beta^K_{1,j}}[-j]  @>{\partial_1}>> R @>>> 0.\\
\end{CD}
$$}

We know that $\phi_0=1_R$ by construction and that $\phi_1$ is a rank
$n$ matrix (over $k$) all of whose entries are in $k$ because $\Y$ is
minimally contained in $I$.  Let $\E_{\bullet}$ denote the mapping
cone on the diagram induced by $\phi$,
$$\E_{\bullet}:=\ 0 \to R \xrightarrow[]{\psi_{n+1}}
R^{\alpha^I_n}\oplus R^n \xrightarrow[]{\psi_n} \cdots
\xrightarrow[]{\psi_2} R^{\alpha^I_1} \oplus R \xrightarrow[]{\psi_1}
R\to 0,$$
where we have used $\alpha^I_j$ to denote the $j$th Betti
number of $R/I$ and have suppressed the graded notation at this step
so that the resolution is more legible. The dual of $\E_{\bullet}$ is
$$\E_{\bullet}^*:=\ 0 \to R \xrightarrow[]{\psi^*_1} R^{\alpha^I_1}
\oplus R \xrightarrow[]{\psi^*_2} R^{\alpha^I_2}\oplus R^n
\xrightarrow[]{\psi^*_3} \cdots \xrightarrow[]{\psi^*_4}
R^{\alpha^I_n}\oplus R \xrightarrow[]{\psi^*_{n+1}} R\to 0,$$
and it
is not difficult to show that $\E_{\bullet}^*$ is a free resolution of
$R/(\Y:I)$. This resolution is never minimal, but we are able to
identify the cause of the non-minimality in the $(n-1)$st, the $n$th,
and the $(n+1)$st terms of $\E_{\bullet}^*$. In fact, the map
$\psi^*_1$ is just multiplication by $1_R$ (actually, $1_{R^*}$) in
the right coordinate, $\psi^*_1(m)=(0,m)$. This implies that the copy
of $R$ constituting $\E_{n+1}$ maps isomorphically onto the copy of
$R$ belonging to $\F_n$ in $\E^*_{n}$, and we may remove both from the
resolution.  So $$\E_{\bullet}' := 0 \rightarrow R^{\alpha^I_1}
\xrightarrow[]{\psi^*_2} R^{\alpha^I_2}\oplus R^n
\xrightarrow[]{\psi^*_3} \cdots \xrightarrow[]{\psi^*_4}
R^{\alpha^I_n}\oplus R \xrightarrow[]{\psi^*_{n+1}} R\to 0,$$
is a
free resolution of $R/(\Y:I)$ where we abuse notation and reuse
$\psi^*_2$ to denote the restriction of $\psi_2^*$ to
$R^{\alpha^I_1}$.

Now for $m\in R^{\alpha^I_1}$, $\psi^*_2(m)=
(\delta^*_2(m),-\phi^*_1(m))$, and as we noted above $\phi_1$ (and
hence also $\phi^*_1$) is a rank $n$ matrix consisting of degree zero
elements.  Thus for each $i$, a copy of $R[-\omega+a_i]$ in $\E'_n$
maps isomorphically onto the copy of $R[-\omega+a_i]$ in $\E'_{n-1}$
(we remember the grading at this step). These pairs may be removed
from $\E_{\bullet}'$, so write $\overline\psi^*_2$ to be the map given
by restriction of $\psi^*_2$ to $\sum_{j}
R^{\beta^I_{1,j}-|j|}[-\omega+j]$, and $\overline\psi^*_3$ to be the
restriction of $\psi^*_3$ to $ \sum_{j} R^{\beta^I_{2,j}}[-\omega+j]$.
Thus
$$\E_{\bullet}'':= 0 \to \sum_{j} R^{\beta^I_{1,j}-|j|}[-\omega+j]
\xrightarrow[]{\overline\psi^*_2} \sum_{j}
R^{\beta^I_{2,j}}[-\omega+j] \xrightarrow[]{\overline\psi^*_3} \cdots
\xrightarrow[]{\psi^*_{n+1}} R\to 0,$$
is a free resolution of $R/I$,
and although it may fail to be minimal, $\overline\psi_2^*$ at least
is rank zero over $k$ (that is, no further cancellation can occur
between $\sum_{j} R^{\beta^I_{1,j}-|j|}[-\omega+j]$ and $\sum_{j}
R^{\beta^I_{2,j}}[-\omega+j]$). We conclude that
$\beta_{n,\omega-j}^{(\Y:I)}=\beta_{1,j}^I-|j|$ as required.

In the case that $I$ fails to minimally contain an $\A$-regular
sequence, this argument needs only a small modification. The cyclic
module $R/(\Y:I)$ can again be resolved using the dual of the mapping
cone on $R/\Y\to R/I$, yielding $$\E_{\bullet}^*:=\ 0 \to R
\xrightarrow[]{\psi^*_1} R^{\alpha^I_1}\oplus R
\xrightarrow[]{\psi^*_2} R^{\alpha^I_2}\oplus R^n
\xrightarrow[]{\psi^*_3} \cdots \xrightarrow[]{\psi^*_4}
R^{\alpha^I_n}\oplus R \xrightarrow[]{\psi^*_{n+1}} R\to 0,$$
and we
can again remove the extra copy of $R$ which constitutes $\E^*_{n+1}$.
The result then follows after noting that there are at most $|j|$
copies of $R[-\omega +j]$ in $\E'_{n-1}=\sum_{j}
R^{\beta^I_{2,j}}[-\omega+j] \oplus \sum_{j}
R^{\beta^K_{1,j}}[-\omega+j]$ which can cancel with copies of
$R[-\omega+j]$ in $\E'_n=\sum_{j} R^{\beta^I_{1,j}}[-\omega+j]$. We
conclude that $\beta_{n,\omega-j}^{(\Y:I)}=\beta_{1,j}^I-t_j$ for
$0\le t_j\le |j|$ as required.
\end{proof}

Proving the main theorem of this section is now easily accomplished.
\begin{theorem}
  Suppose that $L$ is lex plus powers with respect to $\A=\a$,
  $I\subset R$, both share the same Hilbert function, and $I$ contains
  an $\A$-regular sequence.  If the lex plus powers conjecture for
  socles (conjecture \ref{c:socles}) holds, then $\beta^L_{1,j}\ge
  \beta^I_{1,j}$ for all $j$.
\end{theorem}
\begin{proof}
  Let $\X=\{x_1^{a_1}, \dots , x_n^{a_n}\}\subset L$ and let $\Y$ be
  an $\a$-regular sequence in $I$. We know that $(\X:L)$ and $(\Y:I)$
  share the same Hilbert function, the former is $\ap$-lex plus
  powers, and the latter contains an $\ap$-regular sequence (by lemma
  \ref{l:regsequence}). By proposition \ref{p:lasttwo}, $
  \beta^{(\X:L)}_{n,\omega -j}=\beta^L_{1,j}-|j|$, and
  $\beta^{(\Y:I)}_{n,\omega-j}=\beta^I_{1,j}-t_j$. But by hypothesis,
  $\beta^{(\X:L)}_{n,j}\ge \beta^{(\Y:I)}_{n,j}$, and as $|j|\ge t_j$,
  we conclude that $\beta^L_{1,j}\ge \beta^I_{1,j}$ for all $j$ as
  required.
\end{proof}

We conclude by noting that in order to prove conjecture
\ref{c:socles}, it is enough to demonstrate that lex plus powers
ideals have largest socles in a single degree. In particular,
conjecture \ref{c:socles}, and hence LPPH, is equivalent to the
following:

\begin{conjecture}\label{c:soclelast}
  Let $L$ be $\LA$ for some $\A=\a$ and let $\ph$ be the regularity of
  $\H=H(R/L)$. Then $\beta^{L}_{n,\ph+n-1}\ge \beta^I_{n,\ph+n-1}$ for
  any ideal $I\subset R$ containing an $\A$-regular sequence and
  attaining $\H$.
\end{conjecture}

\begin{theorem}
  Conjecture \ref{c:soclelast} and conjecture \ref{c:socles} are
  equivalent.
\end{theorem}
\begin{proof}
  It is obvious that conjecture \ref{c:socles} implies conjecture
  \ref{c:soclelast}.  So suppose that conjecture \ref{c:soclelast}
  holds, $L$ is $\LA$ for some $\A=\a$, $I\subset R$ contains an
  $\A$-regular sequence, and $H(R/L)=H(R/I)=\H$ has regularity $\ph$.
  Now $\beta^L_{n,\ph+n-1}\ge \beta^I_{n,\ph+n-1}$ by hypothesis and
  $\beta^L_{n,\ph+n}=\beta^I_{n,\ph+n}$ because $L$ and $I$ attain the
  same Hilbert function. Thus it remains to show that
  $\beta^L_{n,j}\ge \beta^I_{n,j}$ for all $j \le \ph+n-2$. This is
  easily accomplished. Let $\overline L$ and $\overline I$ be the
  ideals $L+\m{x_1, \dots , x_n}^{\ph}$ and $I+\m{x_1, \dots ,
    x_n}^{\ph}$ respectively. Then $H(R/\overline L)=H(R/\overline I)$
  and $\poh=\ph-1$; by induction on $\rho$ we have that
  $\beta^L_{n,j}=\beta^{\overline L}_{n,j}\ge \beta^{\overline
    I}_{n,j}=\beta^I_{n,j}$ for $j \le \ph + n -2$ as required (where
  we make use of the fact that adding $\m{x_1, \dots , x_n}^{\ph}$ to
  $L$ and $I$ only perturbs the last two rows of their Betti
  diagrams).
\end{proof}

\section*{Acknowledgments} The second author gratefully acknowledges
the financial support received from an NSERC PDF.  The second author
is also grateful for valuable discussions with Susan Cooper and Leslie
Roberts.

\end{document}